\RequirePackage{etex}
\documentclass[final,onefignum,onetabnum]{siamart190516}

\usepackage{lipsum}
\usepackage{amsfonts}
\usepackage{graphicx}
\usepackage{epstopdf}
\ifpdf
  \DeclareGraphicsExtensions{.eps,.pdf,.png,.jpg}
\else
  \DeclareGraphicsExtensions{.eps}
\fi

\newsiamremark{remark}{Remark}
\newsiamremark{hypothesis}{Hypothesis}
\crefname{hypothesis}{Hypothesis}{Hypotheses}
\newsiamthm{claim}{Claim}

\headers{Universal conditional gradient sliding}{Y. Ouyang and T. Squires}

\title{Universal conditional gradient sliding for convex optimization\thanks{Submitted to the editors DATE.
\funding{The authors are partially supported by the Office of Navel Research grant N00014-
20-1-2089.}}}

\author{Yuyuan Ouyang\thanks{School of Mathematical and Statistical Sciences, Clemson University, Clemson, SC
  (\email{yuyuano@clemson.edu}).}
\and Trevor Squires\thanks{School of Mathematical and Statistical Sciences, Clemson University, Clemson, SC
  (\email{tsquire@clemson.edu}).}
  }

\usepackage{amsopn}

\ifpdf
\hypersetup{
  pdftitle={Universal conditional gradient sliding for convex optimization},
  pdfauthor={Y. Ouyang and T. Squires}
}
\fi

\usepackage{amsfonts,amsmath,amssymb}
\usepackage{tsquire}
\usepackage{algorithm, algpseudocode, multirow}
\usepackage{braket}
\usepackage{enumitem}

\usepackage{autonum}
\usepackage{placeins}

\begin{document}

\maketitle

\begin{abstract}
In this paper, we present a first-order projection-free method, namely, the universal conditional gradient sliding (UCGS) method, for solving $\varepsilon$-approximate solutions to convex differentiable optimization problems. For objective functions with H\"{o}lder continuous gradients, we show that UCGS is able to terminate with $\varepsilon$-solutions with at most $\cO((M_\nu D_X^{1+\nu}/\varepsilon)^{2/(1+3\nu)})$ gradient evaluations and $\cO((M_\nu D_X^{1+\nu}/\varepsilon)^{4/(1+3\nu)})$ linear objective optimizations, where $\nu\in (0,1]$ and $M_\nu>0$ are the exponent and constant of the H\"older condition. 
Furthermore, UCGS is able to perform such computations without requiring any specific knowledge of the smoothness information $\nu$ and $M_\nu$. 
In the weakly smooth case when $\nu\in (0,1)$, both complexity results improve the current state-of-the-art $\cO((M_\nu D_X^{1+\nu}/\epsilon)^{1/\nu})$ results \cite{nesterov2018holderCG,ghadimi2019conditional} on first-order projection-free method achieved by the conditional gradient method. Within the class of sliding-type algorithms followed from the work of \cite{lan2016conditional,lan2016gradient}, to the best of our knowledge, this is the first time a sliding-type algorithm is able to improve not only the gradient complexity but also the overall complexity for computing an approximate solution. In the smooth case when $\nu=1$, UCGS matches the state-of-the-art complexity result achieved by the conditional gradient sliding method \cite{lan2016conditional}, but adds more features allowing for practical implementation. 
\end{abstract}

\begin{keywords}
  Convex optimization, first-order method, conditional gradient method, conditional gradient sliding, universal gradient method
\end{keywords}

\begin{AMS}
  90C25, 90C06, 49M37
\end{AMS}

\section{Introduction}
\label{sec:intro}

In this paper, we study first-order projection-free methods for computing $\varepsilon$-approximation solutions to convex optimization problems of form
\begin{equation}\label{problem:CP}
    f^* = \underset{x\in X}{\min}\: f(x).
\end{equation}
Here $X\subset \R^n$ is a high-dimensional compact convex set, $f$ is a convex function, and our goal is to compute an $\varepsilon$-solution $y\in X$ such that $f(y)-f^*\le\varepsilon$. We make the following assumptions concerning the set $X$ and the objective function $f$. Under the Euclidean norm $\|\cdot\|$, we assume that $X$ is compact with diameter
\begin{align}
	\label{eq:D}
	D_X := & \underset{x,y \in X}{\max}\: \norm{x-y} < \infty
\end{align}
and that there exists H\"older exponent $\nu \in (0,1]$ and constant $M_\nu>0$ such that
\begin{equation}
	\label{def:holderCont}
	f(y) \leq f(x) + \ip{\grad f(x),y-x} + \frac{M_\nu}{1+\nu}\norm{x-y}^{1+\nu},\ \forall x,y\in X.
\end{equation}
Our problem of interest covers both smooth ($\nu=1$) and weakly smooth ($\nu\in(0,1)$) convex optimization problems. Specifically, any convex differentiable function whose gradient is $(\nu,M_\nu)-$H\"{o}lder continuous, namely, 
\begin{equation}
    \norm{\grad f(y)- \grad f(x)} \leq M_\nu\norm{y-x}^\nu, \forall x,y \in X
\end{equation}
satisfies \eqref{def:holderCont}.  

Assuming that the projection subproblem $\min_{x\in X}\|x - y\|^2$ can be solved exactly and efficiently for any $y\in\R^n$, first-order projection-based methods for solving the convex optimization problem in \eqref{problem:CP} has already been well studied in the literature. The classical iteration complexity theory  \cite{nemirovski1983problem} has established that the lower complexity bound on the number of gradient evaluations of $\grad f$ is 
\begin{equation}
	\label{eq:nu_complexity}
	\cO(M_\nu D_X^{1+\nu}/\varepsilon)^{2/(1+3\nu)}
\end{equation}
for computing an $\varepsilon$-solution. Note that the above lower complexity bound becomes the widely known $\cO(\sqrt{M_1D_X^2/\varepsilon})$ lower complexity bounds for smooth convex problems (when $\nu=1$). For the smooth case, there have been a large number of literature that developed first-order methods whose performance matches the lower complexity bounds (see, e.g., the books/monographs \cite{bubeck2015convex,beck2017first,nesterov2018lectures,lan2020first} and references with in). 
There also exist several first-order methods in the literature that are able to uniformly achieve the lower complexity bound \eqref{eq:nu_complexity} for any $\nu\in[0,1]$, including for example the fast gradient method (FGM) developed in \cite{nesterov2015universal}, the bundle-level type methods in \cite{lan2015bundle}, and the fast bundle-level method in \cite{chen2019fast} (for the case when $X$ is either a Euclidean norm ball or $\R^n$). Note that the methods in \cite{nesterov2015universal,lan2015bundle,chen2019fast} are universal methods, in the sense that they do not require any knowledge on the values of $\nu$ and $M_\nu$ and are able to achieve the complexity \eqref{eq:nu_complexity} with the best possible $\nu\in[0,1]$ and $M_\nu>0$. Such uniform property is appealing since it allows the methods in \cite{nesterov2015universal,lan2015bundle,chen2019fast} to be applied to convex optimization problems without requiring any smoothness information, i.e., whether the problem is nonsmooth, smooth, or weakly smooth, while still achieve the lower complexity bound with respect to the best smoothness information.


However, it should be noted that we may not always be able to solve the projection subproblem $\min_{x\in X}\|x - y\|^2$ exactly and efficiently. For example, if $X$ is a general polyhedron, then computing the projection with high accuracy would be challenging when the dimension $n$ is large. Recently, there has been  studies on projection-free methods (see, e.g., \cite{jaggi2013revisiting,harchaoui2015conditional,freund2016new}) that replaces the possibly difficult projection subproblem $\min_{x\in X}\|x - y\|^2$ with the easier-to-solve linear objective subproblems $\min_{x\in X} \langle c, x\rangle$. Such methods can be traced back to \cite{frank1956cndg,levitin1966constrained} and are known as the Frank-Wolfe or conditional gradient methods due to their origin. For the smooth case ($\nu=1$) of problem \eqref{problem:CP}, it is shown in \cite{jaggi2013revisiting,harchaoui2015conditional,freund2016new} that the number of gradient evaluation of $\nabla f$ and linear objective optimization subproblems are upper bounded by $\cO(M_1D_X^2/\varepsilon)$. In implementations the linear subproblems can also be solved approximately within certain accuracy while still maintain the same upper complexity bound. Here the number of linear objective optimization subproblems can not be improved; worse-case problem instances that requires solving at least such number of linear objective optimization subproblems has been shown in \cite{jaggi2013revisiting,lan2020first}. For the general case when $\nu\in(0,1]$, universal methods have been developed in \cite{nesterov2018holderCG,ghadimi2019conditional} that compute  $\varepsilon$-solutions with at most $\cO((M_\nu D_X^{1+\nu}/\epsilon)^{1/\nu})$ gradient evaluation of $\grad f$ and linear objective optimization subproblems. 

Focusing on the number of gradient evaluations of $\grad f$ required by the aforementioned projection-free methods, we can observe a significant gap with the lower complexity bound in \eqref{eq:nu_complexity}. For example, the number of gradient evaluations required by the universal methods in \cite{nesterov2018holderCG,ghadimi2019conditional} is upper bounded by $\cO(1/\varepsilon^{3})$ when $\nu=1/3$. This complexity is significantly worse than the lower complexity bound in \eqref{eq:nu_complexity} which is of order $\cO(1/\varepsilon)$ when $\nu=1/3$. For the special smooth case (when $\nu=1$), the number of gradient evaluations required by the methods in \cite{jaggi2013revisiting,harchaoui2015conditional,freund2016new} are upper bounded by $\cO(1/\varepsilon)$ and is significantly worse than the $\cO(1/\sqrt{\varepsilon})$ lower complexity bound in \eqref{eq:nu_complexity}. 

Recently, there has been a breakthrough on closing some of the gap in the gradient evaluations of $\grad f$ between the upper and lower complexity bounds. For the special smooth case (when $\nu=1$), a condition gradient sliding (CGS) method is proposed in \cite{lan2016conditional} that is able to compute an $\varepsilon$-approximate solution of problem \eqref{problem:CP} with $\cO(\sqrt{M_1D_X^2/\varepsilon})$ gradient evaluations of $\nabla f$ and $\cO(M_1D_X^2/\varepsilon)$ linear objective optimization subproblems. Here the number of gradient evaluations required by the CGS method matches that in the lower complexity bound in \eqref{eq:nu_complexity} (with $\nu=1$). Also, the number of linear objective subproblems required by the CGS method also matches the lower bound in \cite{jaggi2013revisiting,lan2020first}. 
Therefore, CGS is the first method that reaches the performance limit of first-order projection-free methods for solving the special smooth case of problem \eqref{problem:CP}. 
It should be noted that the CGS method in \cite{lan2016conditional} requires the knowledge of the Lipschitz continuity constant $M_1$ of the gradient $\grad f$ and does not have a termination criterion for verifying whether it has computed an $\varepsilon$-solution. A backtracking linesearch version of the CGS method is proposed recently in \cite{nazari2020backtracking}, which has the same computational complexity as in \cite{lan2016conditional}, while only requires an initial guess $L_0\le \cO(1)M_1$ and the diameter constant $D_X$ for its computation. It is also able to terminate whenever it verifies the successful computation of an $\varepsilon$-solution. 



None of the above literature on sliding-type algorithms discuss the weakly smooth case when $\nu\in(0,1)$ or the design of universal methods. In this paper, we propose to close the remaining gap in the gradient evaluations of $\grad f$ between its upper complexity bounds in projection-free methods
and the lower complexity bounds in \eqref{eq:nu_complexity}. Specifically, we propose a novel first-order projection-free method, namely the universal conditional gradient sliding (UCGS) method, that is able to compute an $\varepsilon$-solution of the problem \eqref{problem:CP} 
without requiring any knowledge of the smoothness information $(\nu,M_\nu)$. The framework of UCGS is built around that of the fast gradient \cite{nesterov2015universal} and conditional gradient sliding  \cite{lan2016conditional} methods. The contributions of this paper are summarized below.

First, in terms of gradient evaluations of $\grad f$, the total number of evaluations required by the proposed UCGS method for computing an $\varepsilon$-solution is upper bounded uniformly by $\cO(M_\nu D_X^{1+\nu}/\varepsilon)^{2/(1+3\nu)}$ for any $\nu\in(0,1]$. Such bound matches the lower complexity bound in \eqref{eq:nu_complexity}. To the best of our knowledge, this is the first first-order projection-free method that is able to achieve such gradient evaluation complexity bound uniformly for smooth and weakly smooth convex optimization problems.

Second, the total number of linear objective subproblems required by the proposed UCGS method for computing an $\varepsilon$-solution is upper bounded uniformly by $\cO(M_\nu D_X^{1+\nu}/\varepsilon)^{4/(1+3\nu)}$ for any $\nu\in(0,1]$. Comparing with the $\cO((M_\nu D_X^{1+\nu}/\epsilon)^{1/\nu})$ result \cite{nesterov2018holderCG,ghadimi2019conditional} in the literature, the proposed UCGS method has the same complexity when $\nu=1$ and is significantly better for all $\nu\in(0,1)$. For example, when $\nu=1/3$, the UCGS method has significantly better complexity of $\cO(1/\varepsilon^2)$ comparing the $\cO(1/\varepsilon^3)$ result in \cite{nesterov2018holderCG,ghadimi2019conditional}. Within the class of sliding-type algorithms followed from the work of \cite{lan2016conditional,lan2016gradient}, to the best of our knowledge, this is the first time a sliding-type algorithm is able to improve not only the gradient complexity but also the overall complexity for computing an approximate solution.


Third, the proposed UCGS method is able to achieve the aforementioned complexity bounds without any knowledge of the smooth information $(\nu,M_\nu)$ of the objective function. Therefore, it is a universal method that is able to solve weakly smooth and smooth convex optimization problems with the best possible $\nu\in(0,1]$ and $M_\nu>0$.
Note that in the special smooth case when $\nu=1$, the proposed UCGS method can be understood as an extension of the CGS method \cite{lan2016conditional} with add features for practical implementation. In such case, it has the same complexity results as the CGS method \cite{lan2016conditional} and its backtracking linesearch edition \cite{nazari2020backtracking} in terms of both gradient evaluations of $\grad f$ and linear objective subproblems. However, unlike the linesearch edition \cite{nazari2020backtracking}, by incorporating a different backtracking linesearch strategy with a novel parameter choice, UCGS no longer require any information on the continuity constant $M_1$. UCGS also allows that all linear objective optimization subproblems be solved approximately within certain accuracy while maintain the same complexity results.


This paper is organized as follows. In Section \ref{sec:gugAnalysis} we provide a generic description of algorithms for solving problem \eqref{problem:CP}. We discuss the relation of our generic description to the previously mentioned algorithms and perform theoretical complexity analysis.  In particular, we demonstrate the theoretical novelty of our paper on the improvement of both the gradient and linear objective optimization complexities over that of the conditional gradient method. In Section \ref{sec:UCGSAnalysis}, we provide a practical version, namely UCGS, of the generic algorithm and prove that UCGS maintains the same complexity results. We report some preliminary numerical results in Section \ref{sec:numResults} and provide concluding remarks and potential future works in Section \ref{sec:conclusion}.

\section{Conditional gradient sliding method in the H\"older case}
\label{sec:gugAnalysis}
In this section, we propose to analyze the conditional gradient (CG) and conditional gradient sliding (CGS) method in \cite{lan2016conditional} through a generic description, which we call the generic universal gradient (GUG) method. 
GUG can be simply understood as a generic description of the CGS method and serves purely as a tool for our theoretical analysis of CGS in the H\"older case.
We will show that a version of GUG achieves better gradient evaluation complexity than CG for solving problem \eqref{problem:CP}. 
While CGS can already achieve better gradient evaluation than that of CG for problems with Lipschitz continuous gradients, our result covers a more general case of problems with H\"older continuous gradients. 
Moreover, we will also show a novel theoretical result that GUG can also achieve better complexity on linear objective optimizations than that of CG when the H\"older continuouity exponent $\nu\in (0,1)$. Such theoretical result is particularly interesting within the class of sliding-type algorithms followed from the works of \cite{lan2016conditional,lan2016gradient}. To the best of our knowledge, this is the first time a sliding-type algorithm is able to improve not only the gradient complexity but also the overal complexity for computing an approximate solution.

The iterations of GUG is described in Algorithm \ref{alg:GUG}. Let us make a few remarks regarding Algorithm \ref{alg:GUG}. 
First, in both Algorithms \ref{alg:GUG} and in the sequel, we refer to the operations between increments in $t$ as an inner iteration and that of $k$ as an outer iteration. To distinguish inner and outer iteration descriptions, we will use subscripts and superscripts to denote outer and inner iteration indices, respectively. 
Second, the relation \eqref{eq:subproblem} in the Approx-Subproblem procedure can be satisfied through different algorithms and allows both projection-based and projection-free implementations. For example, if we require $\eta_k\equiv 0$, then the Approx-Subproblem procedure solves a projection problem and the iterate $x_k$ computed by the procedure is an optimal solution to the projection problem
\begin{equation}\label{def:phiK}
    \underset{x\in X}{\min}\:{\phi_k(x)}:= \ip{\grad f(z_k),x} + \frac{\beta_k}{2}\norm{x-x_{k-1}}^2.
\end{equation}

\begin{algorithm}
\caption{\label{alg:GUG} Generic universal gradient (GUG) method}
\footnotesize
    \begin{algorithmic}
        \State Start: Choose tolerance $\varepsilon > 0$ and initial iteration $x_0 \in X$.  Set $y_0 = x_0$. 
        \For {$k=1,2,\ldots,N$}
            \begin{align}
                \label{def:zk_gug}
                z_k &= (1-\gamma_k)y_{k-1} + \gamma_k x_{k-1}
                \\
                x_k &= \text{Approx-Subproblem}(\grad f(z_k),x_{k-1},\beta_k, \eta_k)
                \\
                \label{def:yk_gug}
                y_k &= (1-\gamma_k)y_{k-1}+\gamma_kx_k
            \end{align}
        \EndFor
        \State Output $y_N$ as the approximate solution.\\

        \Procedure{$u^+ =$ Approx-Subproblem}{$g,u,\beta,\eta$}
           \State Use any algorithm to compute an approximate solution $u^+$ to problem
           \begin{align}
                \label{def:phi_gug}
                \min_{x\in X}\phi(x):=\ip{g, x} + \frac{\beta}{2}\|x - u\|^2
           \end{align}
           that satisfies
           \begin{equation}\label{eq:subproblem}
                \underset{x\in X}{\max}\:\ip{\grad \phi(u^+), u^+ - x} = \underset{x\in X}{\max}\:\ip{g + \beta(u^+-u), u^+ - x} \leq \eta.          
           \end{equation}
        \EndProcedure
    \end{algorithmic}
\end{algorithm}

\begin{algorithm}
\caption{\label{alg:cgm} Conditional gradient method (CGM) procedure for solving \eqref{eq:subproblem}}
\footnotesize
    \begin{algorithmic}
    \Procedure{$u^+ =$ CGM}{$g,u,\beta$}
    \State Initialize $u^0 = u$. 
                \While {$u^{t-1}$ does not satisfy \eqref{eq:subproblem}}
                    \State Compute $v^t$ such that
                    \begin{equation}\label{def:vt}
                        \underset{x\in X}{\max}\: \ip{g + \beta(u^{t-1} - u), v^t-x}\leq 0
                    \end{equation}
                    \State Set
                    \begin{equation}\label{def:ut}
                        u^t = (1-\alpha^t)u^{t-1} + \alpha^t v^t
                    \end{equation}
                \EndWhile
                \State Output $u^+=u^t$
    \EndProcedure
    \end{algorithmic}
\end{algorithm}

\noindent Consequently, GUG reduces to a version of Nesterov's accelerated gradient method (see, e.g., \cite{nesterov2004introductory}).
For our study, we will focus on a projection-free implementation of the Approx-Subproblem procedure, namely the conditional gradient method (CGM) procedure, described in Algorithm \ref{alg:cgm}. 
Third, if $\beta_k\equiv 0$, then the subproblem \eqref{eq:subproblem} becomes a linear objection optimization and it takes exactly one inner iteration for CGM to compute an optimal solution to this subproblem. Consequently, GUG reduces to CG. Note that by the description of $v^t$ in \eqref{def:vt}, $v^t$ is the optimal solution to the linear subproblem. If instead we allow the right hand side of \eqref{def:vt} to be nonzero, then we can study practical implementation variants of CG that solve the linear objective optimization subproblem approximately (see, e.g., \cite{freund2016new} and the references within). However, we will focus on theoretical analysis in this section; the approximate linear subproblem implementation will be discussed in next section.
Finally, if the parameter $\alpha^t$ in the CGM procedure is chosen as described in \eqref{def:alpha} later, then CGM is exactly CndG in \cite{lan2016conditional}, and GUG reduces to CGS. The key concept behind CGS, which distinguishes it from projection-based methods and CG, is that uses the CGM procedure with multiple inner iterations to compute an approximate solution to the projection problem \eqref{def:phiK}. Instead of solving an optimal solution, CGS runs several inner iterations through the CGM procedure to compute an approximate solution $x_k$ satisfying 
\begin{equation}\label{def:eta}
    \ip{\grad f(z_k) + \beta_k(x_k-x_{k-1}), x_k - x} \leq \eta_k,\ \forall x\in X,
\end{equation}
where $\beta_k>0$.
By doing so, in the special smooth case (when $\nu=1$) of problem \ref{problem:CP} CGS successfully reduces the total number of outer iterations, and hence skipping gradient evaluations, to $\cO(1/\sqrt{\varepsilon})$ from CG's $\cO(1/\varepsilon)$, while not affecting the total number of linear objective optimizations. This feature of skipping gradient evaluations is termed ``sliding'' in \cite{lan2016conditional} (see also \cite{lan2016gradient}). An interesting discovery we will make in this section, is that the use of sliding also reduces the total number of inner iterations, and consequently linear objective optimizations when the H\"older exponent $\nu\in(0,1)$. To the best of our knowledge, such discovery was not made previously in the literature of sliding-type algorithms.

We will now analyze the performance of Algorithm \ref{alg:GUG} under various parameter settings using the CGM procedures in Algorithm \ref{alg:cgm} to solve an approximate solution that satisfies the requirement \eqref{eq:subproblem} in GUG. 
We begin by building a recurrence relation on the outer iterates. Such recurrence provides us a tool for performing complexity analysis on GUG.
\begin{proposition}\label{prop:outerIterGUG}
    Suppose that $\gamma_k \in [0,1]$ for all $k$ in Algorithm \ref{alg:GUG}. We have
    \begin{align}
        \label{eq:f_outer}
        \begin{aligned}
        & f(y_k) - (1-\gamma_k)f(y_{k-1}) - \gamma_k f(x) 
        \\
        \le & \gamma_k\eta_k + \frac{\beta_k\gamma_k}{2}(\norm{x_{k-1}-x}^2 - \norm{x_k-x}^2)
        \\
        & - \frac{\beta_k\gamma_k}{2} \norm{x_k-x_{k-1}}^2+ \frac{M_\nu\gamma_k^{1+\nu}}{1+\nu}\norm{x_k-x_{k-1}}^{1+\nu},\ \forall k\ge 1, x\in X.
        \end{aligned}
    \end{align}
    Specially, if $\nu\in (0,1)$ and $\beta_k>0$ for all $k$, then 
    \begin{align}
        \label{eq:f_outer_xi}
        \begin{aligned}
        & f(y_k) - (1-\gamma_k)f(y_{k-1}) - \gamma_k f(x) 
        \\
        \le & \gamma_k\eta_k + \frac{\beta_k\gamma_k}{2}(\norm{x_{k-1}-x}^2 - \norm{x_k-x}^2) + \xi_k,\ \forall k\ge 1, x\in X,
        \end{aligned}
    \end{align}
    where
    \begin{align}
        \label{def:xik}
        \xi_k :=  \frac{1-\nu}{2(1+\nu)} M_\nu^{\frac{2}{1-\nu}}\left(\frac{\gamma_k}{\beta_k}\right)^{\frac{1+\nu}{1-\nu}}.
    \end{align}
\end{proposition}


\begin{proof}
    From the H\"older condition \eqref{def:holderCont} and the convexity of $f(x)$ we have
    \begin{align}
        & f(y_k) - (1-\gamma_k)f(y_{k-1}) - \gamma_k f(x) 
        \\
        \le & f(z_k) + \ip{\grad f(z_k), y_k - z_k} + \frac{M_\nu}{1+\nu}\|y_k - z_k\|^{1+\nu}
        \\
        & - (1-\gamma_k) (f(z_k) + \ip{\grad f(z_k), y_{k-1} - z_k}) - \gamma_k(f(z_k) + \ip{\grad f(z_k), x - z_k})
        \\
        = & \gamma_k\ip{\grad f(z_k), x_k - x} + \frac{M_\nu\gamma_k^{1+\nu}}{1+\nu}\|x_k - x_{k-1}\|^{1+\nu}.
    \end{align}
    Here the last equality is from the definitions of $z_k$ and $y_k$ in \eqref{def:zk_gug} and \eqref{def:yk_gug} respectively. Noting that $x_k$ is computed from the Approx-Subproblem procedure that satisfies \eqref{eq:subproblem}, we have
    \begin{align}
        & \ip{\grad f(z_k), x_k - x} + \frac{\beta_k}{2}(\|x_k - x_{k-1}\|^2 + \|x_k - x\|^2 - \|x_{k-1}-x\|^2) 
        \\
        = & \ip{\grad f(z_k) + \beta_k(x_k-x_{k-1}), x_k - x} \leq \eta_k,\ \forall x\in X.
    \end{align}
    Summarizing the above two relations we obtain \eqref{eq:f_outer}. By Young's inequality (applied to the product of $(\beta_k\gamma_k/(1+\nu))^{(1+\nu)/2)}\|x_k - x_{k-1}\|^{1+\nu}$ and $M_\nu(\gamma_k/\beta_k)^{(1+\nu)/2}(1+\nu)^{-(1-\nu)/2}$ with exponents $2/(1+\nu)$ and $2/(1-\nu)$ respectively) we conclude the next result\eqref{eq:f_outer_xi} from \eqref{eq:f_outer}.
\end{proof}

In the above proposition there is a recurrence relation concerning weights $(1-\gamma_k)$. The following notation will be used in the sequel for analyzing the complexity of GUG:
    \begin{equation}\label{def:Gammak}
        \Gamma_k = \begin{cases} 1 & k = 1\\ \Gamma_{k-1}(1-\gamma_k) & k > 1.\end{cases}.
    \end{equation}
We will use the following simple lemma for analyzing the sum of recurrent terms.
\begin{lemma}
	\label{lem:tele}
	Suppose that $\{a_k\},\{b_k\}\subset\R$ and $\{\gamma_k\}\subset [0,1]$ are sequences that satisfy $\gamma_1=1$ and
	\begin{align}
		\label{eq:tele_assum}
		a_k \le (1-\gamma_k)a_{k-1} + \gamma_k b_k,\ \forall k\ge 1.
	\end{align}
	Then we have
	\begin{align}
		\label{eq:tele_res}
		a_k\le \Gamma_k\sum_{i=1}^{k}\frac{\gamma_i}{\Gamma_i}b_i,\ \forall k\ge 1,\text{ where }\Gamma_k:=\begin{cases}
			1 & \text{ when }k=1
			\\
			\Gamma_{k-1}(1-\gamma_k) & \text{ when }k>1.
		\end{cases}
	\end{align}
\end{lemma}
\begin{proof}
	Dividing both sides of \eqref{eq:tele_assum} by $\Gamma_k$ we obtain a series of inequalities with telescoping terms concerning sequence $\{a_k/\Gamma_k\}$. Summing up we obtain \eqref{eq:tele_res}.
\end{proof}

We are now ready to derive results on the complexity of CG as a special case of GUG. 
Theorem \ref{thm:cgBound} below is a known complexity result of CG for problems with H\"older continuous gradients  \cite{nesterov2018holderCG,ghadimi2019conditional}.

\begin{theorem}[see also \cite{nesterov2018holderCG,ghadimi2019conditional}]\label{thm:cgBound}
    Suppose that we apply GUG in Algorithm \ref{alg:GUG} (with Algorithm \ref{alg:cgm} to solve Approx-Subproblem) with parameters $\beta_k \equiv 0$, $\eta_k \equiv 0$, $\gamma_k = {2}/{(k+1)}$ and $\alpha^1=1$ in Algorithm \ref{alg:cgm}. To compute an $\varepsilon$-solution to problem \eqref{problem:CP} with H\"older exponent $\nu$ and constant $M_\nu$, GUG requires at most $N_\text{grad}$ gradient evaluations and $N_\text{lin}$ linear objective optimizations, in which
    \begin{equation}\label{eq:boundCG}
        N_\text{lin} = N_\text{grad} = \cO\left(\left( \frac{M_\nu D_X^{1+\nu}}{\varepsilon}\right)^{\frac{1}{\nu}}\right).
    \end{equation}
\end{theorem}
\begin{proof}
    Since $\gamma_k=2/(k+1)$, by \eqref{def:Gammak} we have $\Gamma_k = 2/(k(k+1))$ and hence $\gamma_k/\Gamma_k=k$. Applying Proposition \ref{prop:outerIterGUG} and noting Lemma \ref{lem:tele} with our parameter settings, we have for any $x\in X$ that
    \begin{equation}
        f(y_N) - f(x) \leq \frac{2M_\nu}{N(N+1)(1+\nu)}  \sum_{k=1}^N k\left(\frac{2}{k+1}\right)^{\nu}\|x_k - x_{k-1}\|^{1+\nu}\le \cO\left(\frac{M_\nu D_X^{1+\nu}}{N^\nu}\right).
    \end{equation}
    Thus, in order to obtain an $\varepsilon$-solution, we require 
    at most $N_{grad}$ outer iterations. 
    Moreover, noting that $\alpha^1=1$ and $\beta=0$ in the CGM procedure, comparing \eqref{def:vt} and \eqref{eq:subproblem} we observe that CGM will always terminate after one inner iteration. Therefore, the total number of gradient evaluations and linear optimizations must both be upper bounded by \eqref{eq:boundCG}.
\end{proof}

As pointed in the remarks after the description of Algorithm \ref{alg:GUG}, CG is a special case of GUG with $\beta_k \equiv 0$. Therefore, Theorem \ref{thm:cgBound} above provides a complexity result for the CG algorithm applied to functions with H\"{o}lder continuity exponent $\nu \in (0,1]$. One achieves similar results to Theorem \ref{thm:cgBound} when choosing different $\gamma_k$ (e.g., $\gamma_k=1/k$; see, e.g., \cite{nesterov2018holderCG} for other choices of $\gamma_k$). 
It should also be noted that the choice of $\eta_k\equiv 0$ does not affect the above analysis; indeed, with $\beta=0$ and any $\eta\ge 0$, the CGM procedure will always terminate after one inner iteration. However, as we describe below, if $\beta>0$, the choice of $\eta$ will affect the number of inner iterations performed by the CGM procedure before termination. The proposition below is a known complexity result (see Theorem 2.2(c) in \cite{lan2016conditional}) of CG for solving projection problems. For completeness, we will prove it later in the next section as an immediate consequence of Proposition \ref{prop:innerIterUCGS}.


\begin{proposition}\label{prop:innerIterGUG}
	In the CGM procedure for computing an approximate solution to the projection problem \eqref{def:phi_gug}, if we choose $\alpha^t = 2/(t+1)$, then
    \begin{align}
        \underset{j =0, \dots, t}{\min}\:\underset{x\in X}{\max}\: \ip{\grad \phi(u^{j}), u^{j}-x} \leq \frac{6\beta D_X^2}{t},\ \forall t\ge 1.
    \end{align}
\end{proposition}

Proposition \ref{prop:innerIterGUG} provides insight on the number of inner iterations required by the CGM procedure in Algorithm \ref{alg:cgm} to solve the projection problem \ref{def:phi_gug} approximately. 
For example, if we set $\eta_k \geq 6\beta_k D_X^2$, then the CGM procedure always terminates after exactly one iteration. Noting that $u^0=u$ in CGM, we can observe that GUG reduces to CG not only when $\beta_k\equiv 0$ (as stated previously in the remarks of GUG and after Theorem \ref{thm:cgBound}), but also when $\beta_k>0$ and $\eta_k\ge 6\beta_k D_X^2$. The latter observation is important for our analysis: as described in the following theorem, for problems with H\"older continuous exponent $\nu\in (0,1)$, the latter observation will allow us to perform a simple analysis of CG that is different from the current literature \cite{nesterov2018holderCG,ghadimi2019conditional}. Such simple analysis leads to our interesting discovery that sliding could improve the complexity of linear objective optimizations.

\begin{theorem}[see also \cite{nesterov2018holderCG,ghadimi2019conditional}]\label{thm:cgBound_alt}
    Assume in problem \eqref{problem:CP} that the H\"older exponent $\nu\in(0,1)$. Suppose that we apply GUG in Algorithm \ref{alg:GUG} (with Algorithm \ref{alg:cgm} to solve Approx-Subproblem) with parameters $\beta_k > 0$, $\eta_k = 6\beta_k D_X^2$, and $\alpha^1=1$ in Algorithm \ref{alg:cgm}. Then we have for any $x\in X$ that
    \begin{align}
        \label{eq:cg_cgs_key}
    \begin{aligned}
        & f(y_N) - f(x) 
        \\
        \le & \Gamma_N\sum_{k=1}^{N}\frac{\beta_k\gamma_k}{2\Gamma_k}(12D_X^2 + \norm{x_{k-1}-x}^2 - \norm{x_k-x}^2) + \frac{1}{\Gamma_k}\frac{1-\nu}{2(1+\nu)} M_\nu^{\frac{2}{1-\nu}}\left(\frac{\gamma_k}{\beta_k}\right)^{\frac{1+\nu}{1-\nu}}.
    \end{aligned}
    \end{align}
    Specially, if we set $\beta_k = M_\nu\gamma_k^\nu/D_X^{1-\nu}$ and $\gamma_k = 2/(k+1)$, to compute an $\varepsilon$-solution to problem \eqref{problem:CP} with H\"older exponent $\nu$ and constant $M_\nu$, GUG requires at most $N_\text{grad}$ gradient evaluations and $N_\text{lin}$ linear objective optimizations, in which 
    $
        N_\text{lin} = N_\text{grad} = \cO\left(\left( {M_\nu D_X^{1+\nu}}/{\varepsilon}\right)^{\frac{1}{\nu}}\right).
    $
\end{theorem}
\begin{proof}
    Applying Proposition \ref{prop:outerIterGUG} and noting Lemma \ref{lem:tele}, with our choice of $\eta_k$ we have \eqref{eq:cg_cgs_key}. Consequently,
    \begin{align}
        \label{eq:cg_drawback}
        f(y_N) - f(x)
        \le & \Gamma_N\sum_{k=1}^{N}\frac{7\beta_k\gamma_k D_X^2}{\Gamma_k} + \frac{1}{\Gamma_k}\frac{1-\nu}{2(1+\nu)} M_\nu^{\frac{2}{1-\nu}}\left(\frac{\gamma_k}{\beta_k}\right)^{\frac{1+\nu}{1-\nu}}
        ,\ \forall x\in X.
    \end{align}
    Substituting to \eqref{eq:cg_drawback} the values of $\beta_k$, $\gamma_k$, and noting that $\Gamma_k = 1/k$ from \eqref{def:Gammak}, we have
    \begin{align}
        f(y_N) - f(x)\le \cO\left(\frac{M_\nu D_X^{1+\nu}}{N^2}\right)  \sum_{k=1}^N \frac{k}{(k+1)^\nu}\le \cO\left(\frac{M_\nu D_X^{1+\nu}}{N^\nu}\right).
    \end{align}
    Thus, in order to obtain an $\varepsilon$-solution, we require 
    at most $N_{grad}$ outer iterations. 
    Moreover, noting that $\alpha^1=1$ and $\eta = 6\beta D_X^2$ in the CGM procedure, by Proposition \ref{prop:innerIterGUG} we observe that the CGM procedure will always terminate after one inner iteration. Therefore the total number of linear optimizations is upper bounded by $N_{lin} = N_{grad}$.
\end{proof}

In the above theorem, we observe an imperfection by in the derivation from \eqref{eq:cg_cgs_key} and \eqref{eq:cg_drawback}, although we obtain the same complexity result of CG as in Theorem \ref{thm:cgBound}. Specifically, due to the existence of the dominant term $D_X^2$, we can only simply bound the telescoping difference $(\norm{x_{k-1}-x}^2 - \norm{x_k-x}^2)$ by $D_X^2$ in to obtain \eqref{eq:cg_drawback}. As a consequence, even if we attempt to choose the best $\beta_k = \cO(M_\nu\gamma_k^\nu/D_X^{1-\nu})$ to minimize the right hand side of \eqref{eq:cg_drawback}, the complexity result remains to be $\cO\left(\left( {M_\nu D_X^{1+\nu}}/{\varepsilon}\right)^{\frac{1}{\nu}}\right)$. Noting that the imperfection we observe is due to the choice that $\eta_k = 6\beta_kD_X^2$, we may choose a smaller $\eta_k$ setting to improve the complexity results, as stated in the proposition below.

\begin{theorem}\label{thm:cgsBound}
    Assume in problem \eqref{problem:CP} that the H\"older exponent $\nu\in(0,1)$. Suppose that we apply GUG in Algorithm \ref{alg:GUG} (with Algorithm \ref{alg:cgm} to solve Approx-Subproblem) with parameters 
    \begin{equation}
        \beta_k = \frac{M_\nu k^{\frac{1-3\nu}{2}}}{D_X^{1-\nu}}, \eta_k = \frac{6\beta_k D_X^2}{k}, \text{ and  }\gamma_k = \frac{2}{k+1}.
    \end{equation} 
    To compute an $\varepsilon$-solution to problem \eqref{problem:CP} with H\"older exponent $\nu$ and constant $M_\nu$, GUG requires at most $N_\text{grad}$ gradient evaluations and $N_\text{lin}$ linear objective optimizations, in which
    \begin{equation}
        N_\text{grad} = \cO\left(\left( \frac{M_\nu D_X^{1+\nu}}{\varepsilon}\right)^{\frac{2}{1+3\nu}}\right)\text{ and }
        N_\text{lin} = \cO\left(\left( \frac{M_\nu  D_X^{1+\nu}}{\varepsilon}\right)^{\frac{4}{1+3\nu}}\right).
    \end{equation}
    for any $\nu \in [0,1)$. 
\end{theorem}
\begin{proof}
    Since $\gamma_k = {2}/{(k+1)}$, we have $\Gamma_k = {2}/{(k(k+1))}$ and hence $\gamma_k/\Gamma_k=k$. Applying Proposition \ref{prop:outerIterGUG} and noting Lemma \ref{lem:tele}, with our choice of parameters
    \begin{align}
        & f(y_N) - f(x) 
        \\
        \leq & \frac{2}{N(N+1)} \sum_{k=1}^N 6\beta_k D_X^2 + \frac{k\beta_k}{2}\left(\norm{x_{k-1}-x}^2-\norm{x_k-x}^2\right) + \frac{\xi_k}{\Gamma_k},\ \forall x\in X.
    \end{align}
    Noting that $k\beta_k$ is increasing, we have
    \begin{align}
        & \sum_{k=1}^{N}k\beta_k\left(\norm{x_{k-1}-x}^2-\norm{x_k-x}^2\right) 
        \\
        = & \beta_1\|x_0 - x\|^2 + \sum_{k=1}^{N}((k+1)\beta_{k+1} - k\beta_k)\|x_k - x\|^2 - N\beta_N\|x_N - x\|^2
        \\
        \le & \beta_1 D_X^2 + \sum_{k=1}^{N}((k+1)\beta_{k+1} - k\beta_k)D_X^2 = N\beta_N D_X^2.
    \end{align}
    Combining the above two relations and noting our choice of $\beta_k$ and the description of $\xi_k$ in \eqref{def:xik} we have
    \begin{align}
    \label{eq:balanced}
        f(y_N) - f(x) 
        &\leq \cO\left(\frac{M_\nu D_X^{1+\nu}}{N^2}\right)\left[\sum_{k=1}^{N}k^{\frac{1-3\nu}{2}} + N^{\frac{3-3\nu}{2}} + \sum_{k=1}^{N}\frac{k^{\frac{1+3\nu^2}{2(1-\nu)}}}{(k+1)^{\frac{2\nu}{1-\nu}}} \right] 
        \\
        \le 
        & 
        \cO\left(\frac{M_\nu D_X^{1+\nu}}{N^{\frac{1+3\nu}{2}}}\right),\ \forall x\in X.
    \end{align}
    Thus, to obtain an $\varepsilon$-solution, we need at most $N_{grad}$ outer iterations, or equivalently, at most $N_{grad}$ gradient evaluations. Also, from Proposition \ref{prop:innerIterGUG}, if $\eta_k = {6\beta_kD_X^2}/{k}$, then we will perform at most $k$ inner iterations per outer iteration.  Thus, the total number of inner iterations and consequently linear objective optimizations is upper bounded by
    \begin{equation}
        \sum_{k=1}^{N_\text{grad}} k \le \cO(N_\text{grad}^2) = \cO\left(\left( \frac{M_\nu  D_X^{1+\nu}}{\varepsilon}\right)^{\frac{4}{1+3\nu}}\right).
    \end{equation}
    The proof is now complete.
\end{proof}

Note that by the choice of $\eta_k$ and Proposition \ref{prop:innerIterGUG}, Theorem \ref{thm:cgsBound} provides a complexity result for a version of GUG with the sliding feature for solving problem \eqref{problem:CP}. 
Comparing Theorems \ref{thm:cgBound_alt} and \ref{thm:cgsBound}, the key difference in the proofs is the additional ${1}/{k}$ factor in $\eta_k$ in Theorem \ref{thm:cgsBound}. 
With the additional factor, the three terms at the right hand side of \eqref{eq:balanced} are of the same order with respect to $N$, resolving the imperfection we noticed previously in \eqref{eq:cg_drawback} in the proof of Theorem \ref{thm:cgBound_alt}. In doing so, we achieve the optimal lower complexity bound of gradient evaluations \eqref{eq:nu_complexity} for first-order methods. Interestingly, we can discover that number of linear optimizations required in Theorem \ref{thm:cgsBound} is also significantly reduced comparing with that in Theorem \ref{thm:cgBound_alt}, since $4/(1+3\nu)<1/\nu$
for all $\nu \in (0,1)$. 

It should be noted that we exclude the case $\nu = 1$ case in Theorem \ref{thm:cgsBound} only for convenience of our analysis, since our focus in this section is mainly the theoretical analysis on improving the state-of-the-art complexity bounds \cite{nesterov2018holderCG,ghadimi2019conditional} when $\nu\in (0,1)$. By slightly modifying the proof of Theorem \ref{thm:cgsBound} we can also achieve the same complexity results as the state-of-the-art in \cite{lan2016conditional}. We will include the $\nu=1$ case in the convergence analysis of practical implementation in the next section. 

We conclude this section with several comments regarding the implementation of GUG in Algorithm \ref{alg:GUG}. Note that the sliding result shown in Theorem \ref{thm:cgsBound} requires a parameter choice $\beta_k$ that assumes the knowledge of H\"older exponent $\nu\in (0,1)$ and constant $M_\nu$. Unfortunately, the knowledge of the best $\nu$ and $M_\nu$ for the performance of GUG may not be easily accessible in practice. Furthermore, the proposed Algorithm \ref{alg:GUG} has no termination criterion for verifying whether the current approximate solution $y_k$ is an $\varepsilon$-solution. Lastly, there may exist problem instances in which a solution $v^t$ to the linear subproblem \eqref{def:vt} cannot be computed exactly and instead we can only compute an approximate solution. In the next section, we propose an algorithm called universal conditional gradient sliding (UCGS) that utilizes a backtracking linesearch scheme with an implementable stopping criterion to achieve better practical performance than Algorithm \ref{alg:GUG}. We will also analyze its convergence under an approximate solution to linear subproblem \eqref{def:vt}. 

\section{Practical universal conditional gradient sliding method}
\label{sec:UCGSAnalysis}
In this section, we propose a practical universal conditional gradient sliding (UCGS) method that addresses the above issues of Algorithm \ref{alg:GUG}. The proposed UCGS algorithm is described in Algorithm \ref{alg:UCGS}. 

\begin{algorithm}[t]
	\caption{\label{alg:UCGS} Universal conditional gradient sliding (UCGS) method}
	\footnotesize
	\begin{algorithmic}
		\State
		Start: Choose tolerance $\varepsilon > 0$ and initial iteration $x_0 \in X$.  Set $y_0 = x_0$. 
		\For {$k=1,2,\ldots,$}
		\State Decide $L_k >0$ such that
		\begin{align}
			f(y_k) \leq f(z_k) + \ip{\grad f(z_k),y_k-z_k} + \frac{L_k}{2}\norm{y_k-z_k}^2 + \frac{\varepsilon}{2}\gamma_k\label{def:lk}
		\end{align}
		\State where 
		\begin{align}
			\label{def:gammak} \gamma_k &:= \begin{cases} 1, &k = 1\\ \text{positive solution to } \Gamma_{k-1}(1-\gamma_k)=\frac{L_k\gamma_k^2}{k}, &k > 1 \end{cases}\\
			\label{def:zk} z_k &:= (1-\gamma_k)y_{k-1} + \gamma_k x_{k-1}\\
			\label{def:xk} x_k &:= \text{ACGM}(\grad f(z_k), x_{k-1},\beta_k,\eta_k)\\
			\label{def:yk} y_k &:= (1-\gamma_k)y_{k-1} + \gamma_k x_k\\
			\label{def:Gammakk}\Gamma_k &:= \frac{L_k\gamma_k^2}{k}.
		\end{align}
		
		\State 
		Compute an approximate solution $s_k$ to the problem 
		\begin{align}
		    \label{def:lkproblem}
		    \min_{x\in X}\ell_k(x) := \Gamma_k \sum_{i=1}^k \frac{\gamma_i}{\Gamma_i}\left( f(z_i) + \ip{\grad f(z_i),x-z_i}\right)
		\end{align}
		such that $\ell_k(s_k) - \min_{x\in X}\ell_k(x)\le \varepsilon_k$. Terminate and output $y_k$ as an approximate solution if
		\begin{align}\label{eq:stopCrit}
			f(y_k) - \ell_k(s_k) + \varepsilon_k \leq \varepsilon.
		\end{align}
		\EndFor
		\\

		\Procedure{$u^+ =$ ACGM}{$g,u,\beta, \eta$} 
		\State Goal: Compute $u^+$ such that $\max_{x\in X}\ip{\grad\phi(u^+), u^+-x}\le \eta$, where
		\begin{align}
		    \label{def:phi}
		    \phi(x):=\ip{g, x} + \frac{\beta}{2}\|x - u\|^2.
		\end{align}
		\State Start: Set $u^0 = u$. 
		\For {$t=1,2,\ldots,$}
		\State Compute a $\delta^t$-approximate solution $v^t$ to the problem $\min_{x\in X}\ip{\nabla\phi(u^{t-1}),x}$ such that
		\begin{equation}\label{eq:cndgsubproblem}
			\ip{g + \beta(u^{t-1}-u),v^t} - \min_{x\in X}\ip{g + \beta(u^{t-1}-u), x}  \leq \delta^t.
		\end{equation}
		
		\State Terminate with $u^+:=u^{t-1}$ if
		\begin{align}
		    \label{eq:actual_termi}
		    \ip{g+\beta(u^{t-1}-u),u^{t-1}-v^t} + \delta^t \leq \eta.
		\end{align}
		\State Otherwise, compute
            $
			u^{t}=(1-\alpha^t)u^{t-1}+\alpha^tv^t
			.
			$
		
		\EndFor
		\EndProcedure
	\end{algorithmic}
\end{algorithm}

Let us make a few remarks regarding Algorithm \ref{alg:UCGS}. 
First, the approximate conditional gradient method (ACGM) procedure in Algorithm \ref{alg:UCGS} is a generalization of the CGM procedure (Algorithm \ref{alg:cgm}) discussed in the previous section. Specifically, whenever $\delta^t \equiv 0$, ACGM and CGM are equivalent. Note also that the parameter $\alpha^t$ can be computed through an exact linesearch, namely,
\begin{align}
	\label{def:alpha}
	\alpha^t := \min\: \left\{1,\frac{\ip{g-\beta(u-u^{t-1}),u^{t-1}-v^t}}{\beta\norm{v^t-u^{t-1}}^2}\right\}.
\end{align}
It is easy to observe that the above $\alpha^t$ is the optimal solution to the exact linesearch problem $\min_{\alpha\in[0,1]}\phi((1-\alpha)u^{t-1} + \alpha v^t)$.
Second, if the objective function $f(x)$ in problem \eqref{problem:CP} has Lipschitz continuous gradient (so $\nu=1$) with Lipschitz constant $M_1$, 
then UCGS
can be understood as extension of CGS with added features for practical implementation. The new features include a backtracking linesearch strategy that computes adaptive estimates $L_k$ for the Lipschitz constant $M_1$, the possibility of computing only approximate solutions to linear subproblems, and a termination criterion for verifying whether an approximate solution to problem \eqref{problem:CP} has been computed.
Third,  the choice of $\gamma_k$ and $\Gamma_k$ in \eqref{def:Gammakk} and \eqref{def:gammak} implies that
\begin{equation}\label{eq:GammaRelation}
    \Gamma_k = \begin{cases} 1, &k = 1\\ \Gamma_{k-1}(1-\gamma_k), &k > 1\end{cases}.
\end{equation}
Furthermore, it can be shown that for $k>1$, the solution to $\eqref{def:gammak}$ is given by
\begin{equation}\label{def:solugammak}
    \gamma_k = \frac{2\sqrt{k\Gamma_{k-1}}}{\sqrt{4L_k+k\Gamma_{k-1}}+\sqrt{k\Gamma_{k-1}}}.
\end{equation}
Observe that $\gamma_k \in (0,1)$. Consequently, the recursively described approximate solution $y_k$ is the convex combination of $x_1, \dots, x_k$.  Also the point $z_k$ for gradient evaluation is a convex combination of $x_1, \dots, x_{k-1}$. Such recursive description first appeared in Nesterov's seminal accelerated gradient algorithm (see, e.g., \cite{nesterov2004introductory}) and is also used in the CGS algorithm \cite{lan2016conditional} and the universal gradient algorithms studied in \cite{nesterov2015universal}. However, our choice of $\gamma_k$ is novel and is different from the ones in \cite{nesterov2015universal,lan2016conditional,nesterov2004introductory}. In fact, to our knowledge, none of the settings of $\gamma_k$ in \cite{nesterov2015universal,nesterov2004introductory,lan2016conditional} are suitable for CGS-type algorithms with adaptive $L_k$.  In the only previous work \cite{nazari2020backtracking} that successfully developed a linesearch scheme for CGS, $\gamma_k$ needs to satisfy a more sophisticated cubic equation and $L_k$ needs to be monotone increasing. As we will describe below, such monotonicity restriction on $L_k$ is removed in our proposed UCGS method. 

A few remarks on the practical implementation of Algorithm \ref{alg:UCGS} are also in place.  First, Algorithm \ref{alg:UCGS} proposes that we find $L_k>0$ such that $\eqref{def:lk}$ is satisfied. The condition \eqref{def:lk} originated from the framework of inexact oracle in \cite{devolder2014first} and is also used in \cite{nesterov2015universal}.
We proposed to search for such $L_k$ through a backtracking linesearch strategy. In particular, we initialize with any $L_0 \in \R$ and choose $L_1 = 2^i L_0$ where $i$ is the smallest integer such that \eqref{def:lk} is satisfied. At the start of the $k$-th outer iteration where $k>1$, we set $L_k = L_{k-1}/2$ and assess the validity of $L_k$. If it does not satisfy $\eqref{def:lk}$, we keep backtracking and replacing $L_k$ to $2L_k$ until \eqref{def:lk} is satisfies. Through this backtracking linesearch strategy, we ensure that our choice of $L_k$ is adaptive and that performance is independent of the choice of $L_0$. Previous literature \cite{nazari2020backtracking} on backtracking linesearch strategy of CGS require monotonicity of $L_k$ and may suffer from a poorly chosen $L_0$. Second, our termination criterion is based on \eqref{eq:stopCrit}. We can observe immediately that if the parameter $\varepsilon_k\equiv 0$, i.e., $s_k$ is the exact solution to problem \eqref{def:lkproblem}, then when \eqref{eq:stopCrit} is satisfied, $y_k$ will be $\varepsilon$-approximation solution to problem \eqref{problem:CP}.  To see this, note from \eqref{eq:GammaRelation} that
    \begin{equation}\label{eq:gammaSum}
        \Gamma_k \sum_{i=1}^k \frac{\gamma_i}{\Gamma_i} = 1
    \end{equation}
    and consequently
    $$
        f(y_k) - f^* \le f(y_k) - \min_{x\in X}\ell_k(x) =  f(y_k) - \ell_k(s_k).
    $$
    Such termination criterion also appeared in the previous literature (see, e.g., \cite{nesterov2015universal,nazari2020backtracking}). 
For the case when $\varepsilon_k>0$, we will show later in Theorem \ref{thm:flkDiff} that allowing approximate solution $s_k$ with properly chosen accuracy $\varepsilon_k$ will not affect the complexity results of UCGS.

\FloatBarrier

We present convergence analysis for the UCGS algorithm proposed above, beginning with some results on the inner iteration complexity. The following lemma resembles a combination of the proofs of Theorem 2.2(c) in \cite{lan2016conditional} and Theorem 5.2 in \cite{freund2016new} on the analysis of conditional gradient method with approximate linear objective optimization subproblems for solving projection problems.

\begin{lemma}\label{lem:innerIterUCGS}
    Suppose that $\{\lambda^t\}\in [0,1]$ is any predetermined sequence satisfying $\lambda^1= 1$. In the ACGM procedure, if $\alpha^t$ is chosen such that 
    \begin{align}
        \label{eq:alpha_req}
        \phi(u^t)\le \phi((1-\lambda^t)u^{t-1} + \lambda^tv^t),\ \forall t\ge 1,
    \end{align}
    then we have
    \begin{align}
        \sum_{j=2}^t \frac{\lambda^j}{\Lambda^j}\: \underset{x\in X}{\max}\: \ip{\grad \phi(u^{j-1}),u^{j-1} - x}
        \leq &
        \left[\delta^1 + \sum_{j=2}^t \frac{\lambda^j}{\Lambda^j}\left(\delta^j + \Lambda^{j-1}\sum_{i=1}^{j-1}\frac{\lambda^i}{\Lambda^i}\delta^i\right)\right] 
        \\
        & +
        \frac{\beta D_X^2}{2}\left[1+\sum_{j=2}^t \frac{\lambda^j}{\Lambda^j} \left( \lambda^j + \Lambda^{j-1}\sum_{i=1}^{j-1}\frac{(\lambda^i)^2}{\Lambda^i}\right)\right]
    \end{align}
    for all $t \geq 2$, where
    \begin{align}
        \label{def:Lambda}
        \Lambda^t:=\begin{cases}
            1 & \text{ when }t=1
            \\
            \Lambda^{t-1}(1-\lambda^t) & \text{ when } t>1.
        \end{cases}
    \end{align}
\end{lemma}

\begin{proof}
    Observing that the function $\phi(x)$ in \eqref{def:phi} is a strongly convex function with Lipschitz continuous (with constant $\beta$) gradient, using the assumption \eqref{eq:alpha_req}, and noting the definition of approximate solution $v^t$ in \eqref{eq:cndgsubproblem}, we have
    \begin{align}
    	\label{eq:phiConvex}
    	\begin{aligned}
        &\phi(u^t) - (1-\lambda^t)\phi(u^{t-1}) - \lambda^t(\phi(u^{t-1}) + \ip{\grad \phi(u^{t-1}),x-u^{t-1}})
        \\
        \leq& \phi((1-\lambda^t)u^{t-1} + \lambda^tv^t) - \phi(u^{t-1}) - \lambda^t\ip{\grad\phi(u^{t-1}),x-u^{t-1}}
        \\
        \le & \lambda^t\ip{\grad \phi(u^{t-1}), v^t - u^{t-1}} + \frac{\beta(\lambda^t)^2}{2}\norm{v^t-u^{t-1}}^2 - \lambda^t\ip{\grad \phi(u^{t-1}),x-u^{t-1}}
        \\
        =& \lambda^t\ip{\grad\phi(u^{t-1}),v^t-x} + \frac{\beta(\lambda^t)^2}{2}\norm{v^t-u^{t-1}}^2\\
        \leq& \lambda^t\delta^t +\frac{\beta D_X^2(\lambda^t)^2}{2},\ \forall x\in X, t\ge 1.
	    \end{aligned}
    \end{align}
	Defining $x^*:=\operatorname{argmin}_{x\in X}\phi(x)$, from the above relation we have for any $t\ge 2$ that
	\begin{align}
		\label{eq:tmp1}
		\begin{aligned}
		& \sum_{j=2}^t \frac{\lambda^j}{\Lambda^j} \underset{x\in X}{\max}\:\ip{\grad\phi(u^{j-1}), u^{j-1}-x}
		\\
		\leq& \sum_{j=2}^t \frac{1}{\Lambda^j}[\phi(u^{j-1}) - \phi(x^*)] - \frac{1}{\Lambda^j}[\phi(u^j) - \phi(x^*)] + \frac{\lambda^j}{\Lambda^j}\delta^j + \frac{\beta D_X^2 (\lambda^j)^2}{2\Lambda^j}
		\\
		= & [\phi(u^1) - \phi(x^*)] - \frac{1}{\Lambda^t}[\phi(u^t) - \phi(x^*)] 
		\\
		& + \sum_{j=2}^t \frac{\lambda^j}{\Lambda^j}[\phi(u^{j-1}) - \phi(x^*)] + \frac{\lambda^j}{\Lambda^j}\delta^j + \frac{\beta D_X^2(\lambda^j)^2}{2\Lambda^j}
		\\
		\le & [\phi(u^1) - \phi(x^*)] + \sum_{j=2}^t \frac{\lambda^j}{\Lambda^j}[\phi(u^{j-1}) - \phi(x^*)] + \frac{\lambda^j}{\Lambda^j}\delta^j + \frac{\beta D_X^2(\lambda^j)^2}{2\Lambda^j}.
		\end{aligned}
	\end{align}
	Here in the equality we use the following observations from the definition of $\Lambda^t$ in \eqref{def:Lambda}: $1/\Lambda^1 = 1$ and $1/\Lambda^j = 1/\Lambda^{j-1} + \lambda^j/\Lambda^j$ for all $j\ge 2$, 

    To finish the proof it suffices to bound $\phi(u^{j-1})-\phi(x^*)$ for any $j\ge 2$. Observing that $\phi(x)$ in \eqref{def:phi} is strongly convex and quadratic with
    \begin{equation}
        \frac{\beta}{2}\norm{x-u^{t-1}}^2 = \phi(x) - (\phi(u^{t-1})+\ip{\grad \phi(u^{t-1}), x-u^{t-1}}),\ \forall x\in X, t\ge 1,
    \end{equation}
    we have from \eqref{eq:phiConvex} (with $x = x^*$) that
    \begin{equation}\label{eq:telescopingBound}
        [\phi(u^t)-\phi(x^*)] - (1-\lambda^t)[\phi(u^{t-1})-\phi(x^*)] \leq \lambda^t\delta^t + \frac{\beta D_X^2(\lambda^t)^2}{2} - \frac{\beta\lambda^t}{2}\norm{x^*-u^{t-1}}^2.
    \end{equation}
    Applying Lemma \ref{lem:tele} to the above recurrence relation and ignoring negative terms at the right hand side, we have
    \begin{align}\label{eq:simplePhiBound}
        \phi(u^t) - \phi(x^*) \leq \Lambda^t\sum_{i=1}^t \frac{\lambda^i}{\Lambda^i}\delta^i + \frac{\beta D_X^2(\lambda^i)^2}{2\Lambda^i},\ \forall t\ge 1.
    \end{align}
	We conclude the lemma immediately by applying the above bound to \eqref{eq:tmp1} and rearranging terms.
\end{proof}

The complexity result of the above lemma depends on a predetermined sequence $\{\lambda^t\}$. In the proposition below, we provide a complexity result from an example choice of $\{\lambda^t\}$.

\begin{proposition}\label{prop:innerIterUCGS}
	In the ACGM procedure, at termination we have
	\begin{align}
	    \label{eq:innerIterUCGS_OC}
	    \underset{x\in X}{\max}\: \ip{\grad \phi(u^+), u^+-x} \le \eta.
	\end{align}
	Moreover, if $\delta^t = \sigma\beta D_X^2/t$ for certain $\sigma\ge 0$ and $\alpha^t$ is chosen such that 
	\begin{align}
		\label{eq:alpha_req_2t1}
		\phi(u^t) \le \phi\left(\frac{t-1}{t+1}u^{t-1} + \frac{2}{t+1}v^t\right),
	\end{align}
	then we have for any $t\ge 1$ that
    \begin{align}
        \label{eq:termi_rate}
        \begin{aligned}
        \underset{j=1, \dots, t+1}{\min}\:\ip{\grad \phi(u^{j-1}), u^{j-1}-v^{j}}\le & \underset{j=1, \dots, t+1}{\min}\:\underset{x\in X}{\max}\: \ip{\grad \phi(u^{j-1}), u^{j-1}-x}
        \\
        \le & \frac{6(\sigma+1)\beta D_X^2}{t}.
        \end{aligned}
    \end{align}
    Specially, it takes at most
    \begin{align}
        \label{def:T}
        T := 1 + \left\lceil \frac{(7\sigma+6)\beta D_X^2}{\eta}\right\rceil 
    \end{align}
    iterations for the ACGM procedure to terminate.
\end{proposition}

\begin{proof}
    From the definition of the approximate solution $v^t$ in \eqref{eq:cndgsubproblem}, if the termination criterion in \eqref{eq:actual_termi} of the ACGM procedure is satisfied, then the output $u^+=u^{t-1}$ satisfies
    \begin{align}
        \max_{x\in X}\ip{\grad \phi(u^{t-1}), u^{t-1}-x} &= \max_{x\in X}\ip{\grad \phi(u^{t-1}), u^{t-1}-v^t} + \ip{\grad \phi(u^{t-1}), v^t - x}
        \\
        &\le  (\eta - \delta^t) + \delta^t = \eta.
    \end{align}
    Therefore \eqref{eq:innerIterUCGS_OC} holds. To conclude the proposition it suffices to estimate the rate of convergence of $\underset{x\in X}{\max}\: \ip{\grad \phi(u^{t-1}), u^{t-1}-x}$. To analyze the rate, let us choose $\lambda^t = 2/(t+1)$ and apply Lemma \ref{lem:innerIterUCGS}. Then $\Lambda^t = 2/(t(t+1))$ and
    \begin{align}
        & \sum_{j=2}^t j\cdot\:\underset{x\in X}{\max}\ip{\phi(u^{j-1}),u^{j-1}-x}
        \\
        \leq& \sigma\beta D_X^2\left[1 + \sum_{j=2}^t \left(1 + \frac{2}{j-1}\sum_{i=1}^{j-1} 1\right)\right] + \frac{\beta D_X^2}{2}\left[1+\sum_{j=2}^t\left( \frac{2j}{j+1} + \frac{2}{j-1}\sum_{i=1}^{j-1} \frac{2i}{i+1}\right)\right]
        \\
        < &  \sigma\beta D_X^2(3t-2)+ \frac{\beta D_X^2}{2}(6t-5),\ \forall t\ge 2.
    \end{align}
    Noting that $\sum_{j=2}^t j = (t+2)(t-1)/2$, we have
	\begin{align}
		\underset{j=2, \dots, t}{\min}\:\underset{x\in X}{\max}\: \ip{\grad \phi(u^{j-1}), u^{j-1}-x} \le & \frac{2}{(t+2)(t-1)}\sum_{j=2}^t j\cdot\:\underset{x\in X}{\max}\ip{\phi(u^{j-1}),u^{j-1}-x}
		\\
		< & \frac{6(\sigma+1)\beta D_X^2}{t-1}
		,\ \forall t\ge 2.
	\end{align}
	Using the above result and observing that 
	\begin{align}
		\underset{j=1, \dots, t+1}{\min}\:\ip{\grad \phi(u^{j-1}), u^{j-1}-v^{j}} \le & 
		\underset{j=1, \dots, t+1}{\min}\:\underset{x\in X}{\max}\: \ip{\grad \phi(u^{j-1}), u^{j-1}-x}
		\\
		\le & 
		\underset{j=2, \dots, t+1}{\min}\:\underset{x\in X}{\max}\: \ip{\grad \phi(u^{j-1}), u^{j-1}-x},\ \forall t\ge 1
	\end{align}
	we conclude \eqref{eq:termi_rate}. Moreover, from \eqref{eq:termi_rate} and noting the choice of $\delta^t$, the termination criterion \eqref{eq:actual_termi} holds whenever
    \begin{align}
        \frac{6(\sigma+1)\beta D_X^2}{t-1} + \frac{\sigma \beta D_X^2}{t} \le \eta.
    \end{align}
    Noting the definition of $T$ \eqref{def:T}, the above condition clearly holds for all $t\ge T$.
\end{proof}

In the above proposition, $\sigma\ge 0$ in the definition of $\delta^t$ is a parameter related to the accuracy of approximately solving linear objective optimization subproblems. 
Note that there may also exist other possible choice of $\delta^t$. For example, similar complexity result can be derived by choosing $\delta^t = \sigma \eta$. The benefit of our proposed choice $\delta^t = \sigma \beta D_X^2/t$ from the perspective of practical implementation is that it allows adaptive error of the approximate solution $v^t$ to the linear subproblems and larger error can be admissible when $t$ is small. 

As a side note, recalling that ACGM procedure reduces to CGM procedure in Algorithm \ref{alg:cgm}, we can observe that Proposition \ref{prop:innerIterGUG} in the previous section is a direct consequence of the above result:

\begin{proof}[Proof of Proposition \ref{prop:innerIterGUG}]
	Noting that the CGM procedure described in Algorithm \ref{alg:cgm} is equivalent to the ACGM procedure with $\delta^t\equiv 0$, applying Proposition \ref{prop:innerIterUCGS} above with $\alpha^t = 2/(t+1)$, we conclude the proposition immediately from \eqref{eq:termi_rate}.
\end{proof}

From the above two proofs, it is clear that Proposition \ref{prop:innerIterUCGS} is different from Proposition \ref{prop:innerIterGUG} in the previous section, since it shows us that we can instead compute an approximate solution to \eqref{eq:cndgsubproblem} and proceed with the convergence analysis. We will eventually utilize Proposition \ref{prop:innerIterUCGS} to establish an upper bound on the number of inner iterations that Algorithm \ref{alg:UCGS} requires to compute an $\varepsilon$-solution. We now continue onto the outer iteration analysis, starting with a few results that establish the relation between our computed $L_k$ in the linesearch scheme and the underlying H\"older exponent $\nu$ and constant $M_\nu$ in \eqref{def:holderCont}.
We will use the following lemma that appeared in \cite{nesterov2015universal}. 

\begin{lemma}\label{lem:LkLB}
For any $\delta >0$ and any $L$ such that
\begin{align}
    L \geq \left( \frac{1-\nu}{1+\nu}\cdot\frac{1}{\tau}\right)^{\frac{1-\nu}{1+\nu}} M_\nu ^ \frac{2}{1+\nu},
\end{align} 
where $\nu$ and $M_\nu$ are the H\"older continuity exponent and constant in \eqref{def:holderCont}, we have
\begin{equation}\label{eq:LkUB}
    f(y)\leq f(x) + \ip{\grad f(x), y-x} + \frac{L}{2}\norm{y-x}^2 + \frac{\tau}{2}, \ \forall x,y \in X.
\end{equation}
\end{lemma}
\begin{proof}
    See Lemma 1 of \cite{nesterov2015universal}.
\end{proof}
Note that for $\nu = 1$, the term $\left(\frac{1-\nu}{1+\nu}\right)^\frac{1-\nu}{1+\nu}$ can be handled using a continuity argument $\lim_{\nu\to 1}\left(\frac{1-\nu}{1+\nu}\right)^\frac{1-\nu}{1+\nu} = 1$. We state an immediate corollary of the above lemma below. 
\begin{corollary}\label{coro:LkLB}
    Any $L_k>0$ chosen by Algorithm \ref{alg:UCGS} according to \eqref{def:lk} must also satisfy
    \begin{align}
        L_k \leq 2\left(\frac{1-\nu}{1+\nu}\cdot\frac{1}{\varepsilon \gamma_k}\right)^{\frac{1-\nu}{1+\nu}} M_\nu^{\frac{2}{1+\nu}}
    \end{align}
\end{corollary}
\begin{proof}
    Suppose that $L_k$ does not satisfy \eqref{eq:LkUB}.  Then applying Proposition \ref{prop:flkDiff}  with $\tau = \varepsilon \gamma_k$ implies that $L_k/2$ satisfies \eqref{def:lk}, contradicting the fact that $L_k$ was chosen at step $k$ following the proposed backtracking linesearch implementation (see the remark on practical implementation after the description of Algorithm \ref{alg:UCGS}). 
\end{proof}
The above result is an immediate consequence of the backtracking linesearch strategy we use to find a suitable $L_k$ that satisfies \eqref{def:lk}. 
Based on the above result, we can estimate a bound of $L_k\gamma_k^2$ in the proposition below. Recalling that $\Gamma_k = L_k\gamma_k^2/k$ from \eqref{def:Gammakk} in Algorithm \ref{alg:UCGS}, the following lemma provides also a bound of $\Gamma_k$ that is important for the outer iteration complexity analysis.

\begin{lemma}\label{lem:lkgammak}
    Let $L_k>0$ be chosen by \eqref{def:lk} of Algorithm \ref{alg:UCGS} at step $k$, then
    \begin{align}
        L_k\gamma_k^2 \leq \frac{C_\nu M_\nu^{\frac{2}{1+\nu}}}{k^{\frac{1+3\nu}{1+\nu}}\varepsilon^{\frac{1-\nu}{1+\nu}}} 
    \end{align}
    where 
    \begin{align}\label{eq:c}
        C_\nu := \left( \frac{1+2\nu}{1+3\nu}\right)^{\frac{1+3\nu}{1+\nu}}\left(\frac{1-\nu}{1+\nu}\right)^\frac{1-\nu}{1+\nu}2^{\frac{4+10\nu}{1+\nu}} 
    \end{align}
    is a constant depending only on $\nu$. 
\end{lemma}

\begin{proof}
    The case when $k=1$ is immediate from Corollary \ref{coro:LkLB}. Therefore, throughout the proof we will assume that $k\ge 2$.
    Since we set $\Gamma_k = L_k\gamma_k^2/k$ in Algorithm \ref{alg:UCGS}, we can prove this proposition by bounding $\Gamma_k$. 
    Set $s:= (1+\nu)/(1+3\nu)$.  Since $\nu \in [0,1]$, we have $s\in [{1}/{2},1]$.
    We will study the quantity ${1}/{\Gamma_k^s} - {1}/{\Gamma_{k-1}^s}$, which can be rewritten as
    \begin{align}
        \frac{1}{\Gamma_k^s} - \frac{1}{\Gamma_{k-1}^s} &= \frac{\left( \frac{1}{\Gamma_k^s} - \frac{1}{\Gamma_{k-1}^s}\right)\left( \frac{1}{\Gamma_k^{1-s}} + \frac{1}{\Gamma_{k-1}^{1-s}}\right)}{ \frac{1}{\Gamma_k^{1-s}} + \frac{1}{\Gamma_{k-1}^{1-s}}}\\
        &= \frac{\frac{1}{\Gamma_k}-\frac{1}{\Gamma_{k-1}} - \frac{1}{\Gamma_k}\left(\frac{\Gamma_k}{\Gamma_{k-1}}\right)^s + \frac{1}{\Gamma_{k-1}}\left(\frac{\Gamma_{k-1}}{\Gamma_k}\right)^s}{\frac{1}{\Gamma_k^{1-s}} + \frac{1}{\Gamma_{k-1}^{1-s}}}.
    \end{align}
    Here, noting from the relation of $\Gamma_k$ and $\Gamma_{k-1}$ in \eqref{eq:GammaRelation} that $\Gamma_k \leq \Gamma_{k-1}$ and recalling that $s \in [1/2,1]$, we can make two observations. First, we have $\Gamma_k^{2s-1} \leq \Gamma_{k-1}^{2s-1}$, and hence
    \begin{align}
%
        - \frac{1}{\Gamma_k}\left(\frac{\Gamma_k}{\Gamma_{k-1}}\right)^s + \frac{1}{\Gamma_{k-1}}\left(\frac{\Gamma_{k-1}}{\Gamma_k}\right)^s \geq 0.
    \end{align} 
    Second, we have $\Gamma_{k-1}^{1-s}\le \Gamma_k^{1-s}$, and hence
    \begin{align}
        \frac{1}{\Gamma_k^{1-s}} + \frac{1}{\Gamma_{k-1}^{1-s}}\leq \frac{2}{\Gamma_{k-1}^{1-s}}.
    \end{align}
    Combining the above two observations and recalling that $s=(1+\nu)/(1+3\nu)$ and the relations concerning $\Gamma_k$ and $\Gamma_{k-1}$ in  \eqref{eq:GammaRelation}, we have that
    \begin{align}
        \label{eq:tmp2}
        \frac{1}{\Gamma_k^s} - \frac{1}{\Gamma_{k-1}^s} \geq \frac{\frac{1}{\Gamma_k}-\frac{1}{\Gamma_{k-1}}}{\frac{2}{\Gamma_k^{1-s}}}\nonumber
        = \frac{\frac{\gamma_k}{\Gamma_k}}{\frac{2}{\Gamma_k^{1-s}}}\label{eq:gammaRelation}
        = \frac{\gamma_k}{2}\Gamma_k^{-\frac{1+\nu}{1+3\nu}}\nonumber.
    \end{align}
    We can further bound the last expression in the above relation.  Indeed, recalling that $\Gamma_k = L_k\gamma_k^2/k$ and applying Corollary \ref{coro:LkLB}, we have the inequality
    \begin{align}
        \frac{\gamma_k^2}{k\Gamma_k} = \frac{1}{L_k} \geq \frac{1}{2M_\nu^{\frac{2}{1+\nu}}}\left(\frac{1+\nu}{1-\nu}\cdot\varepsilon \gamma_k\right)^{\frac{1-\nu}{1+\nu}}.
    \end{align}
    In the above recall that we can use a continuity argument for the $\nu=1$ case since $(1-\nu)^{-(1-\nu)}\to 1$ as $\nu\to 1$. Rearranging terms in the above relation, we have
    \begin{align}
        \gamma_k\Gamma_k^{-\frac{1+\nu}{1+3\nu}} \geq \left(\frac{1+\nu}{1-\nu}\right)^\frac{1-\nu}{1+3\nu}\frac{\varepsilon^\frac{1-\nu}{1+3\nu}k^\frac{1+\nu}{1+3\nu}}{2^\frac{1+\nu}{1+3\nu}M_\nu^\frac{2}{1+3\nu}}.
    \end{align}
    Applying the above bound to \eqref{eq:tmp2}, it follows that
    \begin{align}
        \frac{1}{\Gamma_k^s} - \frac{1}{\Gamma_{k-1}^s} \geq \left(\frac{1+\nu}{1-\nu}\right)^\frac{1-\nu}{1+3\nu} \frac{\varepsilon^\frac{1-\nu}{1+3\nu}k^\frac{1+\nu}{1+3\nu}}{2^\frac{2+4\nu}{1+3\nu}M_\nu^\frac{2}{1+3\nu}}.
    \end{align}
    Summing the above from $i = 2$ to $k$ and using the fact that
    \begin{align}
        \sum_{i=2}^k i^{\frac{1+\nu}{1+3\nu}} \ge \int_1^k u^{\frac{1+\nu}{1+3\nu}}du = \frac{1+3\nu}{2+4\nu}\cdot\left(k^{\frac{2+4\nu}{1+3\nu}}-1\right) \ge \frac{1+3\nu}{4+8\nu}k^{\frac{2+4\nu}{1+3\nu}},\ \forall k\ge 2
    \end{align}
    and the definition of $C_\nu$ in \eqref{eq:c}, we obtain
    \begin{align}
        \frac{1}{\Gamma_k^s} \geq \frac{1}{\Gamma_k^s} - \frac{1}{\Gamma_{1}^s} 
        &
        \geq 
        \left(\frac{1+\nu}{1-\nu}\right)^\frac{1-\nu}{1+3\nu}\frac{\varepsilon^\frac{1-\nu}{1+3\nu}}{2^\frac{4+10\nu}{1+3\nu}M_\nu^\frac{2}{1+3\nu}}\frac{1+3\nu}{1+2\nu}k^{\frac{2+4\nu}{1+3\nu}}
    \end{align}
    Recalling that $s=(1+\nu)/(1+3\nu)$ and $\Gamma_k = L_k\gamma_k^2/k$, we conclude the proposition immediately from the above result.
    %
\end{proof}

It should be noted that the technique utilized in Lemma \ref{lem:lkgammak} is similar to that of the proof surrounding equation (4.4) in \cite{nesterov2015universal}. However, note that the choice of parameter $\gamma_k$ in UCGS is different from the one in \cite{nesterov2015universal}.  Therefore, the proof in \cite{nesterov2015universal} needs to be adapted to the above proof. With the help of Lemma \ref{lem:lkgammak}, we are now ready to prove our primary convergence properties on the proposed UCGS algorithm. We start with the following proposition that resembles the outer iteration analysis in Proposition \ref{prop:outerIterGUG} of the previous section.

\begin{proposition}\label{prop:flkDiff}
    Suppose that the parameters in Algorithm \ref{alg:UCGS} satisfy $\beta_k \geq L_k \gamma_k$ for all $k$. Then for any $x\in X$, 
    \begin{align}
        f(y_k) - \ell_k(x) \leq \frac{\varepsilon}{2} + \Gamma_k \sum_{i=1}^k \frac{\gamma_i\beta_i}{2\Gamma_i}\left(\norm{x-x_{i-1}}^2 - \norm{x-x_i}^2\right)^2 + \Gamma_k \sum_{i=1}^k \frac{\gamma_i\eta_i}{\Gamma_i}.
    \end{align} 
\end{proposition}

\begin{proof}
    Fix any $x\in X$. From the definitions of $\ell_k(x)$ and $y_k$ in \eqref{def:lk} and \eqref{def:yk} respectively, we have
    \begin{align}
        \frac{1}{\Gamma_k}\ell_k(x) =& \sum_{i=1}^k \frac{1}{\Gamma_i}\left( \gamma_i f(z_i) + \gamma_i \ip{\grad f(z_i),x-x_i} + \ip{\grad f(z_i),\gamma_i(x_i-z_i)}\right)\\
        =& \sum_{i=1}^k \frac{1}{\Gamma_i}\left( f(z_i) + \ip{\grad f(z_i), y_i-z_i}\right)+ \frac{\gamma_i}{\Gamma_i} \ip{\grad f(z_i),x-x_i}
        \\
        &\qquad- \frac{1-\gamma_i}{\Gamma_i}(f(z_i) + \ip{\grad f(z_i), y_{i-1}-z_i})
    \end{align}
    We will now bound three terms in the above relation.  First, by convexity of $f$, 
    \begin{align}
        -(f(z_i) + \ip{\grad f(z_i),y_{i-1}-z_i}) \ge -f(y_{i-1}).
    \end{align}
    Second, by our choice of $L_k$ in \eqref{def:lk} and the definitions of $y_k$ and $z_k$ in \eqref{def:yk} and \eqref{def:zk} respectively, we have
    \begin{align}
        f(z_i) + \ip{\grad f(z_i),y_i-z_i} \ge & f(y_i) - \frac{L_i}{2}\norm{y_i-z_i}^2 - \frac{\varepsilon}{2}\gamma_i 
        \\
        = & f(y_i) - \frac{L_i\gamma_i^2}{2}\norm{x_i-x_{i-1}}^2 - \frac{\varepsilon}{2}\gamma_i.
    \end{align}
    Lastly, using the result \eqref{eq:innerIterUCGS_OC} in Lemma \ref{lem:innerIterUCGS} and noting the definition of $\phi(x)$ in \eqref{def:phi}, we obtain the following result during the termination of the ACGM procedure in computing $x_i$:
    \begin{align}
         \ip{\grad f(z_i), x-x_i} \geq&\: \beta_i\ip{x_i-x_{i-1},x_i-x} - \eta_i\\
        =& -\frac{\beta_i}{2}\left(\norm{x-x_{i-1}}^2 - \norm{x_i-x_{i-1}}^2 - \norm{x-x_i}^2\right) - \eta_i\\
        \geq& -\frac{\beta_i}{2}\left(\norm{x-x_{i-1}}^2 - \norm{x-x_i}^2\right) -\eta_i + \frac{L_i^2}{2}\norm{x_i-x_{i-1}}^2.
    \end{align}
    In the last inequality above we use our assumption that $\beta_k\geq L_k\gamma_k$ for all $k$. Based on the above three observations and rearranging terms we obtain that
    \begin{align}
        \frac{1}{\Gamma_k}\ell_k(x) \geq& \sum_{i=1}^k \frac{1}{\Gamma_i}\left( f(y_i) - (1-\gamma_i)f(y_{i-1}) - \frac{\gamma_i\beta_i}{2}(\norm{x-x_{i-1}}^2-\norm{x-x_i}^2)\right)
        \\
        &\qquad -\frac{1}{\Gamma_i}\left(\frac{\varepsilon}{2}\gamma_i + \gamma_i\eta_i\right)
        \\
        = & 
        \frac{f(y_k)}{\Gamma_k} -
        \sum_{i=1}^k \frac{\gamma_i\beta_i}{2\Gamma_i}(\norm{x-x_{i-1}}^2 - \norm{x-x_i}^2) -\sum_{i=1}^k \frac{\gamma_i\eta_i}{\Gamma_i} - \frac{\varepsilon}{2\Gamma_k}.
    \end{align}
    Here in the last equality we use the relations \eqref{eq:GammaRelation} and \eqref{eq:gammaSum} and the fact that $\gamma_1=1$ in its definition \eqref{def:gammak}. We conclude the result by multiplying by $\Gamma_k$ and rearranging  terms.
\end{proof}

With the help of Propositions \ref{prop:innerIterUCGS}, Lemma \ref{lem:lkgammak}, and Proposition \ref{prop:flkDiff}, we are ready to present the complexity results of UCGS in the following theorem.

\begin{theorem}\label{thm:flkDiff}
    Suppose that we apply UCGS described in Algorithm \ref{alg:UCGS} with parameters
    \begin{equation}\label{eq:param}
        \beta_k = L_k\gamma_k,\ \eta_k = \frac{L_k\gamma_kD_X^2}{k},\text{ and }\varepsilon_k = \frac{\sigma L_k\gamma_kD_X^2}{2},
    \end{equation}
    and $\alpha^t$ in \eqref{def:alpha} and $\delta^t = \sigma\beta D_X^2/t$ in the ACGM procedure, where $\sigma\ge 0$ is a parameter related to the accuracy of approximately solving linear objective optimization subproblems.
    Then 
    Algorithm \ref{alg:UCGS} terminates with an $\varepsilon$-solution after at most $N_{\text{grad}}$ gradient evaluations and $N_{\text{lin}}$ linear objective optimizations, where
    \begin{align}\label{eq:evaluations}
    \begin{aligned}
        N_{\text{grad}} :=&\left\lceil16\left(\frac{(3+\sigma)^{\frac{1+\nu}{2}}M_\nu D_X^{1+\nu}}{\varepsilon}\right)^\frac{2}{1+3\nu}\right\rceil
        \text{ and } 
        \\
        N_{\text{lin}} :=& \left\lceil\left(\frac{7}{2}\sigma + 3\right)N_{\text{grad}}^2 + \left(\frac{7}{2}\sigma + 6\right)N_{\text{grad}}\right\rceil.
    \end{aligned} 
    \end{align}
\end{theorem}

\begin{proof}
    From the definition of $s_k$ in Algorithm \ref{alg:UCGS}, we have that if the termination criterion of UCGS in \eqref{eq:stopCrit} holds, then $y_k$ is an $\varepsilon$-solution to problem \eqref{problem:CP}. Let us evaluate the number of gradient evaluations, or equivalently, the number of outer iterations of UCGS in order to compute an $\varepsilon$-solution $y_k$.
    Applying Proposition \ref{prop:flkDiff} with our choice of parameters we have
    \begin{align}
        \label{eq:flkStep1}
        \begin{aligned}
        f(y_k) - \ell_k(x) 
        =&\: \frac{\varepsilon}{2} + \frac{\Gamma_k}{2}\sum_{i=1}^k i(\norm{x-x_{i-1}}^2 - \norm{x-x_i}^2) + \Gamma_kkD_X^2,\ \forall x\in X.
        \end{aligned}
    \end{align}
    The second term above can be further simplified by noting from the compactness of $X$ and the definition of the diameter $D_X$ in \eqref{eq:D}. Indeed, we have
    \begin{align}
    \label{eq:telescope}
    \begin{aligned}
        & \sum_{i=1}^k i\left(\norm{x-x_{i-1}}^2 - \norm{x-x_i}^2\right)
        \\
        =& \norm{x-x_0}^2 + \sum_{i=2}^k (i-(i-1))\norm{x-x_{i-1}}^2 - k\norm{x-x_k}^2 
        \\
        \leq&\: D_X^2 + \sum_{i=2}^k D_X^2 =\: kD_X^2,\ \forall x\in X.
    \end{aligned}
    \end{align}
    Thus, we may continue by applying \eqref{eq:telescope} to \eqref{eq:flkStep1} with $x = s_k$ to conclude that
    \begin{align}
        f(y_k) - \ell_k(s_k) \leq \frac{\varepsilon}{2} + \frac{\Gamma_k}{2}k D_X^2  + \Gamma_kkD_X^2 = \frac{\varepsilon}{2} + \frac{3L_k\gamma_k^2}{2} D_X^2.
    \end{align}
    Here we make use of the description of $\Gamma_k$ in \eqref{def:Gammakk} for the last equality. In view of the above result and the value of parameter $\varepsilon_k$ in \eqref{eq:param}, $y_k$ satisfies the termination criterion \eqref{eq:stopCrit} of UCGS and hence becomes an $\varepsilon$-solution whenever $k$ satisfies the relation 
    $(3+\sigma)L_k\gamma_k^2 D_X^2 \leq \varepsilon$. Applying Lemma \ref{lem:lkgammak}, it follows that such relation holds whenever
    \begin{align}
        \frac{(3+\sigma)C_\nu D_X^2  M_\nu^{\frac{2}{1+\nu}}}{k^{\frac{1+3\nu}{1+\nu}}\varepsilon^{\frac{1-\nu}{1+\nu}}} \leq \varepsilon,\text{ i.e., }k\ge C_\nu^\frac{1+\nu}{1+3\nu}\left(\frac{(3+\sigma)^{\frac{1+\nu}{2}}M_\nu D_X^{1+\nu}}{\varepsilon}\right)^\frac{2}{1+3\nu}.
    \end{align}
    Noting that $C_\nu$ defined in \eqref{eq:c} is a constant that depends only on $\nu\in[0,1]$ and observing that $C_\nu^\frac{1+\nu}{1+3\nu} \leq 16$ for all $\nu \in [0,1]$, we conclude that whenever $k \geq N_{\text{grad}}$, $y_k$ is an $\varepsilon$-solution. Therefore, UCGS requires at most $N_{\text{grad}}$ gradient evaluations of $\grad f$ to compute an $\varepsilon$-solution.
    
    It suffices to compute the number of linear objective optimizations that UCGS requires for computing an $\varepsilon$-solution. This is equivalent to estimating the total number of inner iterations that UCGS requires.
    Let us estimate the maximal number of inner iterations required before the termination criterion \eqref{eq:actual_termi} is satisfied.
    Recall from the remark after \eqref{def:alpha} that $\alpha^t$ is the best linesearch parameter and hence satisfies assumption \eqref{eq:alpha_req_2t1} in Proposition \ref{prop:innerIterUCGS}. Applying Proposition \ref{prop:innerIterUCGS} and noting the definition of approximate solution $v^t$ in \eqref{eq:cndgsubproblem}, we have that
    the maximal number of linear objective optimizations performed at the $k$-th call to the ACGM procedure is at most 
    \begin{equation}
        T_k := 1 + \left\lceil \frac{(7\sigma+6)\beta_k D_X^2}{\eta_k}\right\rceil = 1 + \lceil (7\sigma+6)k \rceil.
    \end{equation}
    Adding one linear objective optimization problem in \eqref{def:lkproblem} in each other iteration concerning the termination criterion of UCGS, we conclude that the total number of linear objective optimizations for UCGS to compute an $\varepsilon$-solution is bounded above by
    \begin{align}
        \sum_{k=1}^{N_{\text{grad}}} (T_k + 1) \le  \sum_{k=1}^{N_{\text{grad}}} 3 + (7\sigma+6)k = \left(\frac{7}{2}\sigma + 3\right)N_{\text{grad}}^2 + \left(\frac{7}{2}\sigma + 6\right)N_{\text{grad}}.
    \end{align} 
\end{proof} 

We conclude this section with a few remarks on the above complexity results of UCGS. 
First, we note that UCGS is similar to FGM in \cite{nesterov2015universal} in the sense that the number of gradient evaluations generalizes the accelearated gradient descent method in \cite{nesterov2004introductory}. From Theorem \ref{thm:flkDiff}, number of gradient evaluations required by UCGS to compute an approximate solution is $\cO((M_{\nu}D_X^{1+\nu}/\varepsilon)^{\frac{2}{1+3\nu}})$.  In the smooth case when $\nu = 1$, this becomes $\cO(\sqrt{M_{1}D_X^2/\varepsilon})$ which matches the complexities of gradient evaluations in n\cite{nesterov2015universal,nesterov2004introductory}. 
Second, unlike FGM that requires exact solutions to projection subproblems, we have a bound on the number of linear objective optimizations required to solve the projection subproblem \eqref{def:xk}. From this perspective, UCGS is a generalization of CGS in \cite{lan2016conditional} as a universal method that covers not only the smooth case (when $\nu=1$) but also the weakly smooth case (when $\nu\in(0,1))$, without requiring any knowledge of H\"older exponent $\nu$ and constant $M_\nu$ . Indeed, when $\nu = 1$ our complexity on the number of linear objective optimizations is on the order of $\cO(M_1 D_X^2/\varepsilon)$, which matches the that of CGS in \cite{lan2016conditional}. 
Third, the number of linear objective optimizations and gradient evaluations when CG is applied to \eqref{problem:CP} was shown to be $\cO\left(\left(M_\nu D^{1+\nu}/\varepsilon\right)^{\frac{1}{\nu}}\right)$ in \cite{nesterov2018holderCG,ghadimi2019conditional}. In view of Theorem \ref{thm:flkDiff}, UCGS benefits from sliding and only requires $\cO\left(\left(M_\nu D^{1+\nu}/\varepsilon\right)^{\frac{2}{1+3\nu}}\right)$ gradient evaluations and $\cO\left(\left(M_\nu D^{1+\nu}/\varepsilon\right)^{\frac{4}{1+3\nu}}\right)$ linear objective optimizations which are both improvements over the results in \cite{nesterov2018holderCG,ghadimi2019conditional} whenever $\nu \in (0,1]$. 
Fourth, our proposed UCGS method is not only a universal method generalization of the CGS in \cite{lan2016conditional}. Indeed, there are more features added for practical implementation: it has an implementable exit criterion and allows for an approximate solution to \eqref{eq:cndgsubproblem}. Note that such added features of the UCGS does not affect its theoretical complexity.
Finally, we use the same accuracy constant $\sigma$ for approximately solving the linear subproblems in setting the parameters $\varepsilon_k$ and $\delta^t$. It is easy to change the proof if we use different accuracy constants for $\varepsilon_k$ and $\delta^t$.

\section{Numerical results}
\label{sec:numResults}

Our goal in this section is to present preliminary results from our numerical experiments. We will compare the performance of our proposed UCGS algorithm with that of the CG  method in \cite{nesterov2018holderCG} in two numerical experiments described below. The experiments are performed using MATLAB R2018b.

In the first experiment, we consider the problem
\begin{align}
    \underset{x \in \conv(V)}{\min}\: f(x):= \norm{Ax-b}_2
\end{align}
with $V = \{v_1,\dots,v_p\} \subseteq \R^n$, $\conv(V):= \{x \in \R^n: \exists \lambda \in \Delta_p \text{ s.t. } x=\sum_{j = 1}^p \lambda_i v_i\}$, and $\Delta_p := \{ \lambda \in \R^p: \sum_{i=1}^p \lambda_i = 1, \lambda_i \geq 0\}$ is the standard simplex. In this experiment, we generated vectors $v_i$ uniformly in $[0,1]^n$.  The matrix $A\in \R^{m \times n}$ is a Gaussian randomly generated sparse matrix with density $d$. For this experiment, we fix the number of vectors in the set $V$ to be $p=500$ and set $m = 2n$. The linear objective optimization subproblem is a linear program over the standard simplex and can be computed easily.

For our second experiment, we solve the problem 
\begin{align}
    \underset{X \in \text{Spe}_n}{\min}\: f(X):= \sum_{i=1}^m \norm{X-A_i}_2
\end{align}
where $\text{Spe}_n:= \{X \in \R^{n\times n}: \tr(X) = 1, X\succeq 0\}$ and $A_i \in \text{Spe}_n$ for each $i =1,\dots, m$.  The matrices $A_i$ are obtained by randomly generating an $n\times n$ matrix whose entries follow uniform $[0,1]$ distributions and then projecting it into $\text{Spe}_n$. 
The linear objective optimization problem over $\text{Spe}_n$ is equivalent to a smallest eigenvalue problem, which will be solved by MATLAB's eigs() function. Note that a solution to the smallest eigenvalue problem will not be exact, and therefore we benefit from being able to solve the linear subproblems approximately.

In our experiments, UCGS will terminate whenever an $\varepsilon$-solution with tolerance $\varepsilon = 10^{-3}$ is computed. We will terminate CG if its computational time exceeds twice the amount that UCGS spent before termination. Note that both models in the experiments have nonsmooth objective functions, but are still differentiable at many feasible points.  Therefore, they may benefit from a universal method for $\nu\in (0,1]$. 

\begin{table}[h]
    \begin{center}
        \begin{tabular}{@{}|cc|cccc|ccc|@{}}
            \hline
            & & \multicolumn{4}{|c|}{\textbf{UCGS}} &\multicolumn{3}{|c|}{\textbf{CG}} \\
            $n$ & $d$ &\textbf{GE} &\textbf{LO} & \textbf{Time} & \textbf{Error} & \textbf{Iter} &\textbf{Time} &\textbf{Error}\\ \hline
            $2500$ & $0.2$ & $66$ & $2690$ & $6.71$ & $9.945e-4$ & $572$ & $13.42$ & $9.7086e1$\\

            $2500$ & $0.4$ & $60$ & $3679$ & $9.08$ & $9.976e-4$ & $524$ & $18.17$ & $1.404e2$\\
            $2500$ & $0.6$ & $62$ & $245$ & $2.64$ & $9.678e-4$ & $146$ & $5.29$ & $5.598e2$\\
            $2500$ & $0.8$ & $57$ & $3176$ & $8.45$ & $9.768e-4$ & $399$ & $16.93$ & $2.400e2$\\
            \hline
            $5000$ & $0.2$ & $71$ & $286$ & $7.13$ & $9.882e-4$ & $178$ & $14.32$ & $6.037e2$\\
            $5000$ & $0.4$ & $42$ & $52$ & $4.89$ & $9.585e-4$ & $84$ & $9.81$ & $1.689e3$\\
            $5000$ & $0.6$ & $68$ & $4564$ & $36.14$ & $9.727e-4$ & $483$ & $72.40$ & $3.527e2$\\
            $5000$ & $0.8$ & $67$ & $419$ & $12.91$ & $9.815e-4$ & $161$ & $25.94$ & $1.165e3$\\
            \hline
            $10000$ & $0.2$ & $85$ & $12269$ & $150.51$ & $9.96e-4$ & $915$ & $301.21$ & $2.449e2$\\
            $10000$ & $0.4$ & $69$ & $12614$ & $157.39$ & $9.916e-4$ & $636$ & $315.27$ & $4.734e2$\\
            $10000$ & $0.6$ & $70$ & $16063$ & $205.87$ & $9.821e-4$ & $653$ & $412.14$ & $5.423e2$\\
            $10000$ & $0.8$ & $69$ & $12707$ & $180.65$ & $9.862e-4$ & $473$ & $361.73$ & $8.162e2$\\
            \hline
        \end{tabular}
    \end{center}
    \caption{    \label{table:conv}
Minimizing over a convex hull. Here, we report the gradient evaluations (outer iterations) and linear objective optimization (inner iterations) for UCGS as well as the error that it terminated with. For CG, we allow it to run for twice the amount of time that UCGS took. We then report the number of iterations and whether terminating objective value was better than that of UCGS. }
\end{table}

\begin{table}[h]
    \begin{center}
        \begin{tabular}{@{}|cc|cccc|ccc|@{}}
            \hline
            & & \multicolumn{4}{|c|}{\textbf{UCGS}} &\multicolumn{3}{|c|}{\textbf{CG}} \\
            $n$ & $m$ &\textbf{GE} &\textbf{LO} & \textbf{Time} & \textbf{Error} & \textbf{Iter} &\textbf{Time} &\textbf{Error}\\ \hline
            $50$ & $50$ & $1354$ & $8493$ & $9.87$ & $9.992e-4$ & $6908$ & $19.74$ & $6.073e-3$\\
            $50$ & $100$ & $1767$ & $11138$ & $13.09$ & $9.994e-4$ & $7038$ & $26.19$ & $1.172e-2$\\
            $50$ & $200$ & $2425$ & $15173$ & $25.39$ & $9.995e-4$ & $8273$ & $50.79$ & $2.271e-2$\\
            \hline
            $100$ & $50$ & $1836$ & $13056$ & $159.61$ & $9.980e-4$ & $11648$ & $319.25$ & $3.225e-3$\\
            $100$ & $100$ & $2347$ & $16816$ & $216.59$ & $9.990e-4$ & $13372$ & $433.20$ & $5.634e-3$\\
            $100$ & $200$ & $3296$ & $23836$ & $310.16$ & $9.984e-4$ & $16053$ & $620.36$ & $9.892e-3$\\
            \hline
            $200$ & $50$ & $1722$ & $33673$ & $470.71$ & $9.989e-4$ & $15966$ & $941.43$ & $3.308e-3$\\
            $200$ & $100$ & $2314$ & $46323$ & $730.69$ & $9.994e-4$ & $17033$ & $1461.42$ & $6.870e-3$\\
            $200$ & $200$ & $3154$ & $64511$ & $1086.42$ & $9.992e-4$ & $19762$ & $2172.85$ & $1.015e-2$\\
            \hline
        \end{tabular}
    \end{center}
    \caption{    \label{table:spec}
Minimizing over standard spectrahedron. Here, we report the gradient evaluations (outer iterations) and linear objective optimization (inner iterations) for UCGS as well as the error that it terminated with. For CG, we allow it to run for twice the amount of time that UCGS took. We then report the number of iterations and whether terminating objective value was better than that of UCGS. }

\end{table}

The results from the numerical experiments are documented in Tables \ref{table:conv} and \ref{table:spec}. Columns 1 indicates the sizes $n$\footnote{Note that the length of the vectors are $n$ and $n^2$ in the first and second experiments respectively.} whereas the second column represents either the density of $A$ or the value of $m$ for experiments 1 and 2 respectively. Columns 3 and 4 denote the number of outer iterations, i.e. gradient evaluations (GE), and inner iterations, i.e. linear objective optimization (LO), respectively that UCGS performed before terminating with the desired tolerance.  Columns 5 and 6 present the time (in seconds) used and error upon termination of UCGS. For CG, we report the total number of iterations (Iter) performed, the computational time (in seconds) required and the final error in Columns 7, 8, and 9.  Note that if the time of CG is twice that of UCGS, then the error is not expected to be below our specified tolerance.

Let us make a few comments regarding the results in Tables \ref{table:conv} and \ref{table:spec}.  For the convex hull experiment in Table \ref{table:conv}, we see that the excessive number of gradient evaluations of CG prevents it from being competitive. The gradient of our objective function requires a matrix multiplication of increasingly dense matrices.  As these densities tend to $1$, the gradient evaluations become more computationally expensive, and CG cannot report as good of a solution as UCGS with even in twice the allotted time, because it requires much more gradient evaluations to compute an approximate solution. We also note the necessity of a projection-free algorithm for this feasible set since the projection onto the convex hull requires the solving of a quadratic program. For any moderately sized $n$, this quadratic program is computationally infeasible to solve. For example, one iteration of FGM in \cite{nesterov2015universal} applied to the problem instance with $n=2500$ and $d=0.2$ takes at least $20$ seconds, which is three times as long as UCGS took to converge.

The second experiment over the standard spectrahedron removes the previous difficulty of computing the gradient.  In experiment 2, the cost of the gradient evaluation is almost negligible.  However, we still see in Table \ref{table:spec} that UCGS outperforms CG.  In this case, the superior linear objective optimization complexity of UCGS can be seen by noting that CG performs 1 linear objective optimization per iteration.  Thus, even with a comparable amount of linear objective optimizations, CG can still not match the complexity of UCGS. This directly highlights the differences in the linear objective optimization complexity mentions previously. We also observe the effectiveness of the implementable stopping criterion which enabled us to terminate when an $\varepsilon$-solution was achieved.

\section{Concluding remarks}
\label{sec:conclusion}
In this paper, we present a novel projection-free method, namely the universal conditional gradient sliding (UCGS) method, for convex differentiable optimization with H\"older continuous gradients. We show that UCGS is a generalization of other conditional gradient type methods in terms of gradient evaluations and linear objective optimizations and at the same time has a more practical implementation by requiring less problem dependent parameters. 
Specifically, for an objective function whose gradient is H\"{o}lder continuous with exponent $\nu\in (0,1]$ and constant $M_{\nu}>0$, we prove that UCGS is able to terminate and output an $\varepsilon$-solution with at most $\cO((M_\nu D_X^{1+\nu}/\varepsilon)^{2/(1+3\nu)})$ gradient evaluations and $\cO((M_\nu D_X^{1+\nu}/\varepsilon)^{4/(1+3\nu)})$ linear objective optimizations. Moreover, it is able to perform the computation without requiring any specific knowledge of the H\"older exponent $\nu$ and constant $M_\nu$. UCGS improves the state-of-the-art complexity results (achieved by the conditional gradient method \cite{nesterov2018holderCG,ghadimi2019conditional}) of first-order projection-free methods when $\nu\in (0,1)$. It also matches the state-of-the-art complexity results when $\nu=1$ (achieved by the conditional gradient sliding method \cite{lan2016conditional}) and adds more features allowing for practical implementation. 

The results of this paper can be further generalized. First, if there exists a prox-function $x\mapsto V(u,x)$ defined over the feasible set $X$ that is strongly convex and has Lipschitz continuous gradients with respect to a general norm $\|\cdot\|$ (see, e.g., \cite{ghadimi2012optimal,dang2015convergence} for discussions on such prox-functions), then it can be shown that the UCGS can be generalized to a non-Euclidean version that still achieves $\cO(1/\varepsilon^{2/(1+3\nu)})$ gradient evaluations and $\cO(1/\varepsilon^{4/(1+3\nu)})$ linear objective optimizations. The constants in these complexities will depend on the smoothness information $\nu$ and $M_\nu$ with respect to the general norm $\|\cdot\|$, the diameter of $X$ in terms of $\max_{x,y\in X}V(x,y)$, and the strongly convex constant of the prox-function and the Lipschitz constant of its gradient. 
Second, while we focus on convex differentiable optimization, the results in this paper can also be extended to non-differentiable cases. Indeed, note that a convex function that is Lipschitz continuous with constant $2M_0$ also satisfies the H\"older condition we describe in \eqref{def:holderCont} with exponent $\nu=0$ and constant $M_0$. Also, our proof throughout this paper can be easily extended by replacing gradients $\grad f$ to subgradients $f'$. Consequently, Theorem \ref{thm:flkDiff} shows that UCGS computes an $\varepsilon$-solution with $\cO((M_0 D_X/\varepsilon)^2)$ subgradient evaluations and $\cO((M_0 D_X/\varepsilon)^4)$ linear objective optimizations. 

We conclude this paper by discussing some potential future work. First, similar to the conditional gradient sliding method \cite{lan2016conditional},  
our proposed UCGS requires some information on the diameter $D_X$ of the feasible set $X$. It is interesting to study whether this diameter can be estimated adaptively by backtracking search, similar to the linesearch strategy we use to estimate the smoothness information $L_k$. Second, while we mention in the above generalization of UCGS that our results can be extended to convex non-differentiable optimization, if we focus solely on non-differentiable cases, there exists interesting complexity results in the literature. 
Specifically, in \cite{thekumparampil2020projection} it is shown that in the non-differentiable case it is possible to achieve  $\cO((M_0 D_X/\varepsilon)^2)$ complexity for both subgradient evaluations and linear objective optimizations whenever the objective function is Lipschitz continuous. 
The method in \cite{thekumparampil2020projection} can only be applied to nonsmooth cases when $\nu=0$, and is unknown whether it can be generalized to a universal method that uniformly computes approximate solutions for nonsmooth, weakly smooth, and smooth convex optimization problems. It is a potential future work to study whether the technique in \cite{thekumparampil2020projection} can improve our developed linear objective optimization complexity of the H\"older continuous gradient case when $\nu \in (0,1]$, and whether a universal method for all cases $\nu \in [0,1]$ can be developed.



\bibliographystyle{siamplain}
\bibliography{main}
\end{document}